\documentclass[final,3p,times]{elsarticle}
\usepackage{amsmath,amssymb}
\usepackage{bm}
\usepackage{multirow}
\usepackage{subfigure}% Support for small, `sub' figures and tables
\usepackage{lineno}
\biboptions{sort&compress}
\usepackage{tikz}
\usetikzlibrary{positioning,decorations.pathmorphing}
\tikzset{zigzag/.style={decorate,decoration=zigzag},>=stealth}
\usepackage{makecell}
\usepackage{url}

\begin{document}
	
	\begin{frontmatter}
		
	\title{A Spectral Method for Depth-Separated Solution of a Wavenumber Integration Model in Horizontally Stratified Fluid Acoustic Waveguides}
		
	\author[1]{Houwang Tu}
	\ead{tuhouwang@nudt.edu.cn}
	\author[2]{Yongxian Wang\corref{cor1}}
	\cortext[cor1]{Corresponding author: \url{yxwang@nudt.edu.cn}}
	\author[2]{Wei Liu}
	\ead{liuwei@nudt.edu.cn}
	\author[2]{Shuqing Ma}
	\author[1]{Xiaodong Wang}

	\address[1]{College of Computer, National University of Defense Technology, Changsha, 410073, China}
	\address[2]{College of Meteorology and Oceanography, National University of Defense Technology, Changsha, 410073, China}
		
		\begin{abstract}
			The wavenumber integration model is considered to be the most accurate algorithm for arbitrary horizontally stratified media in computational ocean acoustics. In contrast to the normal mode approach, it considers not only the discrete wavenumber spectrum but also the continuous spectrum components, eliminating errors in the model approximation for horizontally stratified media. Traditionally, analytical and semianalytical methods have been used to solve the depth-separated wave equation in the wavenumber integration method, and numerical solutions have generally focused on the finite difference method and the finite element method. In this paper, an algorithm for solving the depth equation using the Chebyshev--Tau spectral method combined with a domain decomposition strategy is proposed, and a numerical program named WISpec is developed accordingly. The proposed algorithm can simulate not only the sound field excited by a point source but also the sound field excited by an infinite line source. The key idea of the algorithm is to first discretize the depth equations for each layer via the Chebyshev--Tau spectral method and then solve the equations for each layer simultaneously by incorporating boundary and interface conditions. Several representative numerical experiments are devised to test the accuracy and speed of WISpec. The high consistency of the results of different software programs running under the same configuration proves that the numerical algorithm proposed in this paper is accurate, reliable and numerically stable.
		\end{abstract}
		
		%%Graphical abstract
		%\begin{graphicalabstract}
		%\includegraphics{grabs}
		%\end{graphicalabstract}
		
		\begin{keyword}
			Spectral method\sep
			Wavenumber integraion\sep
			Underwater sound propagation\sep
			Computational ocean acoustics
		\end{keyword}
		
	\end{frontmatter}
	
	\newpage
	\section{Introduction}
	The wavenumber integration method, which is one of the most well-established methods in ocean acoustics, is basically an implementation of the integral transform technique for horizontally stratified media \cite{Jensen2011}. This method, which does not make any approximations to the Helmholtz equation, completely avoids approximation error and is considered the most accurate method for simulating sound propagation in horizontally stratified media. The normal mode model is often compared with the wavenumber integration method because the mathematical basis of both is the same; the difference is that the evaluation of the integral adopts different strategies. For the normal mode model, complex contour integration is used to reduce the integral to a sum of residues, whereas in the wavenumber integration method, the integrals are evaluated directly by numerical quadrature \cite{Jensen2011,Etter2018}. The wavenumber spectrum of a general waveguide includes both discrete and continuous components. The discrete wavenumber spectrum can be represented as a sum of normal modes, while the continuous spectrum should be described by an integral in the wavenumber domain. In other words, the normal modes contain only a limited number of discrete wavenumbers that contribute greatly to the sound field while ignoring the continuous spectrum that may have a great error on the sound field, especially the near field. For the case where the horizontal wavenumbers are near the branch cut, the normal mode model may fail to find the root, thus reducing the accuracy of the sound field. Therefore, the wavenumber integration method is generally considered to be more accurate than the normal mode model.
	
	The idea of wavenumber integration for horizontally stratified media was first introduced into ocean acoustics by Pekeris in 1948 \cite{Pekeris1948}. He used simple two- and three-layer structures to model sound propagation in horizontally stratified media. Later, Ewing, Jardetzky and Press used this method to study seismic propagation in waveguides with few layers \cite{EJP1957}. The wavenumber integration technique performs a series of integral transformations on the Helmholtz equation to simplify the original partial differential equation into a series of ordinary differential equations of depth coordinates. These equations are then solved analytically in each layer in such a way that the amplitudes are initially undetermined; these undetermined amplitudes are determined by matching boundary conditions at the interfaces, and finally, the corresponding sound field is solved by computing the inverse integral transform. For the initially proposed ocean environment with few layers, it is easy to solve the linear equations analytically by expressing the boundary conditions in terms of undetermined sound field amplitudes. However, for more complicated ocean environments, the undetermined coefficients method is not applicable, and numerical methods are usually employed.
	
	The earliest algorithm for simulating depth-dependent sound fields is the propagator matrix approach (PMA) proposed by Thomson \cite{Thomson1950} and Haskell \cite{Haskell1953}. The advantage of the PMA is that it is recursive and thus requires only a small amount of memory, but the disadvantage is that it requires a very time-consuming correction scheme to ensure numerical stability. Furthermore, the PMA is not well suited to problems in which the field must be solved at multiple receiver depths \cite{Jensen2011}. Kennett reviewed the PMA \cite{Kennett2009} and proposed the invariant embedding approach (IEA) \cite{Kennett1974}. The advantages of the IEA are its intrinsic numerical stability, the simplicity of the recurrence algorithms and its direct suitability for reflectivity modeling. However, the IEA is not well suited for solving the global problem of interest in ocean acoustics \cite{Kennett1979}. At present, the most widely used method for solving depth equations is the direct global matrix (DGM) approach proposed by Schmidt \cite{Schmidt1985}. In the DGM approach, the sound field of each layer is expressed as the superposition of the sound field generated by the sound source and an undetermined sound field satisfying the homogeneous depth equation, and the relationship between the sound fields of each layer is controlled by the continuity conditions at the interfaces. Then, the depth equations in local layers are assembled into a DGM, and after boundary conditions are imposed, the sound field in all layers can be obtained simultaneously by solving the global linear equations \cite{SAFARI,OASES}. The main advantage of the DGM approach is its unconditional stability, which is obtained at no additional computational cost and yields efficient numerical solutions to the depth-separated wave equations in all layers simultaneously \cite{Schmidt1986}. The problem size and workload of the DGM are proportional to the number of layers. When the acoustic parameters vary greatly with depth or the frequency of the sound source is very high, a denser configuration of layers is required. In this case, the global matrix generated by the DGM approach is larger, and the computational speed is slower \cite{Jensen1998}.
	
	Among the methods for numerically solving differential equations, in addition to the widely used finite difference and finite element methods, spectral methods provide a niche but efficient new tool. Spectral methods offer high accuracy \cite{Orszag1972,Gottlieb1977,Canuto1988} and fast convergence speed \cite{Guoby1998,Boyd2001,Canuto2006,Jshen2011} and have been rapidly developed in acoustics \cite{Sabatini2019,Tuhw2021c}, especially computational ocean acoustics. In recent years, new algorithms for normal modes \cite{Dzieciuch1993,Tuhw2020a,Tuhw2020b,Tuhw2021d,NM-CT,rimLG,MultiLC}, coupled modes \cite{Tuhw2022a,Tuhw2022b,Tuhw2022c}, adiabatic modes \cite{Tuhw2023a} and parabolic equation models \cite{Tuhw2021a,Tuhw2021b,SMPE} based on spectral methods have been successively devised. In the present paper, the Chebyshev--Tau spectral method is used to numerically solve the depth-separated wave equation. In the model designed in this paper, the Chebyshev--Tau spectral method does not physically discretize the ocean environment in the vertical direction; that is, it does not use piecewise linear approximation to address the ocean environmental parameters, so there is no error due to physical discretization. In addition, the algorithm has no factors that make the solution divergent, so it is favorably stable. A corresponding numerical program is developed for the algorithm. Several classic numerical experiments verify the accuracy and illustrate the capability of the algorithm and program devised in this article.
	
\section{Mathematical Modeling}
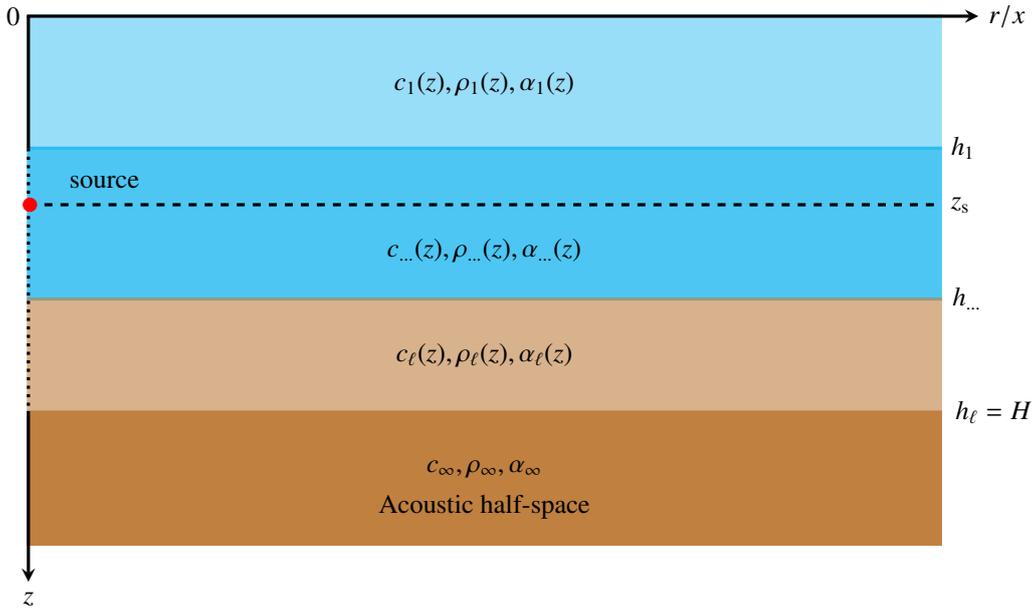
\begin{figure}[htbp]
	\centering
	\begin{tikzpicture}[node distance=2cm]
		\filldraw[very thick,cyan,opacity=0.4] (2,0)--(14,0)--(14,-1.75)--(2,-1.75)--cycle;
		\filldraw[very thick,cyan,opacity=0.7] (2,-1.75)--(14,-1.75)--(14,-3.75)--(2,-3.75)--cycle;
		\filldraw[very thick,brown,opacity=0.6] (2,-3.75)--(14,-3.75)--(14,-5.25)--(2,-5.25)--cycle;
		\filldraw[very thick,brown] (2,-5.25)--(14,-5.25)--(14,-7)--(2,-7)--cycle;		
		\node at (1.8,0){$0$};
		\draw[very thick, ->](1.98,0)--(14.5,0) node[right]{$r/x$};
		\draw[very thick, ->](2,-5.25)--(2,-7.5) node[below]{$z$};
		\draw[very thick](2,0)--(2,-1.75);
		\draw[dotted, very thick](2,-1.75)--(2,-5.25);		
		\node at (14.3,-1.75){$h_1$};
		\node at (14.35,-3.75){$h_{\dots}$};
		\node at (14.7,-5.25){$h_\ell= H$};
		\draw[dashed, very thick](2,-2.5)--(14,-2.5) node[right]{$z_\mathrm{s}$};
		\filldraw [red] (2.02,-2.5) circle [radius=2.5pt];
		\node at (3,-2.2){$\text{source}$};
		\node at (8,-0.9){$c_1(z),\rho_1(z),\alpha_1(z)$};
		\node at (8,-3.1){$c_{\dots}(z),\rho_{\dots}(z),\alpha_{\dots}(z)$};
		\node at (8,-4.5){$c_{\ell}(z),\rho_{\ell}(z),\alpha_{\ell}(z)$};
		\node at (8,-6){$c_{\infty},\rho_{\infty},\alpha_{\infty}$};
		\node at (8,-6.5){Acoustic half-space};
	\end{tikzpicture}
	\caption{Schematic of an arbitrary horizontally stratified ocean environment.}
	\label{Figure1}
\end{figure}
For a horizontally stratified ocean environment, the interfaces at different depths are all parallel planes, and the layer properties are functions only of depth and are independent of the azimuthal angle, as shown in Fig.~\ref{Figure1}. For this two-dimensional range-independent problem, the Helmholtz equation can be written in the following form \cite{Jensen2011}:
\begin{equation}
	\label{eq.1}
	\left[\rho(z)\nabla\cdot\left(\frac{1}{\rho(z)}\nabla\right)+k^{2}(z)\right] \psi(r,z)=F(r,z)
\end{equation}
where $\psi(r,z)$ denotes the displacement potential, $F(r,z)$ is the body force, $k$ is the wavenumber, $k=2\pi f/c(1+\mathrm{i}\eta\alpha)$, $\eta=(40\pi \log_{10}{\mathrm{e}})^{-1}$, $f$ is the frequency of the source, and $c$ and $\alpha$ are the acoustic speed and attenuation of the medium, respectively. The derivations of the sound fields for point and line sources discussed below are based on the Helmholtz equation.

\subsection{Integral transformation for point source problems}
Waveguides excited by point sound sources are usually solved in cylindrical coordinates. The sound field is related only to the depth and horizontal range away from the sound source, so in a cylindrical coordinate system, we let the $z$-axis pass vertically downward through the sound source, and the $r$-axis is parallel to the sea surface, as illustrated in Fig.~\ref{Figure1}. The Helmholtz equation (Eq.~\eqref{eq.1}) in the cylindrical coordinate system is taken in the following form:
\begin{equation}
	\label{eq.2}
	\left[\frac{1}{r} \frac{\partial}{\partial r} \left(r \frac{\partial}{\partial r}\right)+\rho(z) \frac{\partial}{\partial z}\left(\frac{1}{\rho(z)} \frac{\partial }{\partial z}\right)+k^{2}(z)\right] \psi(r, z)=- \frac{\delta(r)\delta\left(z-z_{\mathrm{s}}\right)}{2 \pi r} 
\end{equation}
where $z_\mathrm{s}$ is the depth of the sound source. We consider using the following Hankel transform pairs for the above equation:

\begin{subequations}
	\label{eq.3}
	\begin{align}
		\label{eq.3a}
		\psi(r,z)&=\int_{0}^{\infty} \Psi\left(k_{r}, z\right) J_{0}\left(k_{r} r\right) k_{r} \mathrm{d} k_{r} \\
		\Psi\left(k_{r}, z\right)&=\int_{0}^{\infty} \psi(r, z) J_{0}\left(k_{r} r\right) r \mathrm{d} r
	\end{align}
\end{subequations}
Specifically, the following operation is applied to Eq.~\eqref{eq.2}:
\[
\int_{0}^{\infty} (\cdot) J_{0}\left(k_{r} r\right) r \mathrm{d} r
\]
Therefore, we can easily obtain the following depth-separated wave equation:
\begin{equation}
	\label{eq.4}
	\left[\rho(z)\frac{\mathrm{d}}{\mathrm{d} z}\left(\frac{1}{\rho(z)}\frac{\mathrm{d}}{\mathrm{d} z}\right)+\left(k^{2}-k_{r}^{2}\right)\right] \Psi\left(k_{r}, z\right)=-\frac{\delta\left(z-z_{\mathrm{s}}\right)}{2 \pi}
\end{equation}
This equation is an ordinary differential equation in depth and can be solved numerically or analytically. Conventionally, the solution strategy for the Green function $\Psi(k_r,z)$ is to first physically discretize the ocean environment in the depth direction \cite{Thomson1950,Haskell1953,Kennett2009,Schmidt1986}. The ocean environment is divided into sufficiently thin layers, and the acoustic parameters of each layer are regarded as depth-independent constants, which introduces error. In this paper, we introduce a Chebyshev--Tau spectral method to numerically solve the depth-separated wave equation, which is a high-precision numerical method that does not involve physical discretization. Once $\Psi(k_r,z)$ has been determined at various discrete wavenumbers for the selected receiver depths, Eq.~\eqref{eq.3a} can be evaluated to yield the total displacement potential $\psi(r,z)$ at any depth and range.

\subsection{Integral transformation for line source problems}
An infinitely long line sound source is often used to verify the accuracy of models in computational ocean acoustics. We also consider the solution of this common model. The line source problem is usually introduced in a Cartesian coordinate system, still letting the $z$-axis pass through the sound source with vertical downward orientation; the $x$-axis is parallel to the sea surface, and the sound source penetrates the $xoz$-plane perpendicularly to infinity. The main structure is still as shown in Fig.~\ref{Figure1}, except that the $r$-axis is now the $x$-axis. Therefore, the Helmholtz equation of the line source in the Cartesian coordinate system can be written in the following form \cite{Jensen2011}:
\begin{equation}
	\label{eq.5}
	\left[\frac{\partial^{2}}{\partial x^{2}}+\rho(z) \frac{\partial}{\partial z}\left(\frac{1}{\rho(z)} \frac{\partial }{\partial z}\right)+k^{2}(z)\right] \psi(x, z)=-\delta(x) \delta\left(z-z_{\mathrm{s}}\right)
\end{equation}
We apply the following Fourier transform pairs to Eq.~\eqref{eq.5}:

\begin{subequations}
	\label{eq.6}
	\begin{align}
		\label{eq.6a}
		\psi(x, z) &=\int_{-\infty}^{\infty} \Psi\left(k_{x}, z\right) \mathrm{e}^{\mathrm{i} k_{x} x} \mathrm{~d} k_{x} \\
		\Psi\left(k_{x}, z\right) &=\frac{1}{2 \pi} \int_{-\infty}^{\infty} \psi(x, z) \mathrm{e}^{-\mathrm{i} k_{x} x} \mathrm{d} x
	\end{align}
\end{subequations}
Specifically, the following operator is applied to the above formula:
\[
\frac{1}{2 \pi} \int_{-\infty}^{\infty} (\cdot) \mathrm{e}^{-\mathrm{i} k_{x} x} \mathrm{d} x
\]
The following depth-separated wave equation is thus obtained:
\begin{equation}
	\label{eq.7}
	\left[\rho(z)\frac{\mathrm{d}}{\mathrm{d} z}\left(\frac{1}{\rho(z)}\frac{\mathrm{d}}{\mathrm{d} z}\right)+\left(k^{2}-k_{x}^{2}\right)\right] \Psi\left(k_{x}, z\right)=- \frac{\delta\left(z-z_{\mathrm{s}}\right)}{2 \pi}
\end{equation}
Solving Eq.~\eqref{eq.7} yields the depth-dependent Green function $\Psi(k_x,z)$. After obtaining $\Psi(k_x,z)$, the total sound field can be synthesized by Eq.~\eqref{eq.6a}, as discussed for the point source.

A comparison of Eqs.~\eqref{eq.4} and \eqref{eq.7} indicates that the depth-separated wave equations for the point source and line source have exactly the same form, except that $r$ is replaced by $x$ and $k_r$ is replaced by $k_x$. Therefore, the Green function of the depth-separated wave equation can be used as the integral kernel function not only for the point source but also for the line source. We take only Eq.~\eqref{eq.4} as an example for the solution of the depth-separated wave equation above.

\subsection{Interface conditions and boundary conditions}
In the ocean environment shown in Fig.~\ref{Figure1}, the interfaces ($\{h_l\}_{l=1}^{\ell-1}$) with discontinuous environmental parameters in seawater need to satisfy the interface conditions. According to the relationship between sound pressure and displacement potential, the sound pressure and normal particle velocity must be continuous, yielding:

\begin{subequations}
	\label{eq.8}
	\begin{gather}
		\rho(h_l^+)\Psi(k_r,h_l^+)=\rho(h_l^-)\Psi(k_r,h_l^-),\quad l=1,2,\cdots,\ell-1\\
		\label{eq.8b}
		\left.\frac{\mathrm{d} \Psi}{\mathrm{d} 	z}\right|_{z=h_l^+}=\left.\frac{\mathrm{d} \Psi}{\mathrm{d} z}\right|_{z=h_l^-},\quad l=1,2,\cdots,\ell-1
	\end{gather}
\end{subequations}
where the superscripts $-$ and $+$ indicate the interfaces from above and below, respectively.

The definite solution of Eq.~\eqref{eq.4} includes the boundary conditions at the sea surface ($z=0$) and the seabed ($z=H$). Considering the large difference in impedance between seawater and air, the sea surface is usually taken as a perfectly reflected boundary, that is, a pressure-release boundary:
\begin{equation}
	\label{eq.9}
	\psi(r,0)=0 \Longleftrightarrow\Psi(k_r,0) = 0
\end{equation}
For the lower boundary condition, pressure-release seabed and rigid seabed are usually used:

\begin{subequations}
	\label{eq.10}	
	\begin{gather}
		\psi(r,H)=0 \Longleftrightarrow\Psi(k_r,H) = 0\\
		\left.\frac{\mathrm{d} \psi_{\ell}(r,z)}{\mathrm{d} z}\right|_{z=H}=0 \Longleftrightarrow \left.\frac{\mathrm{d} \Psi_{\ell}(k_r,z)}{\mathrm{d} z}\right|_{z=H} = 0
	\end{gather}
\end{subequations}
In addition, when modeling the ocean environment, the acoustic half-space boundary shown in Fig.~\ref{Figure1} is also typically used in practice. The lower boundary condition that needs to be satisfied on the truncated bottom of the acoustic half-space is:
\begin{equation}
	\label{eq.11}
	\rho_\infty \left.\frac{\mathrm{d} \Psi_\ell}{\mathrm{d} z}\right|_{z=H} + \left.\rho_\ell(H)\sqrt{k_r^2-k_\infty^2} \Psi_\ell\right|_{z=H} = 0
\end{equation}

Note that the inhomogeneous term on the right-hand side of Eq.~\eqref{eq.4} contains the source term $\delta(z-z_\mathrm{s})$, so the singularity necessitates special treatment of the equation at the depth of the sound source. We add a virtual interface at the depth of the sound source. 

Due to the singularity of the sound source, the normal particle velocity cannot be constrained using the continuity condition Eq.~\eqref{eq.8b}. A natural idea is to integrate both sides of Eq.~\eqref{eq.4} in a very small neighborhood of $z_\mathrm{s}$ so that $\delta(z-z_\mathrm{s})$ can be eliminated.
\begin{equation}
	\label{eq.12}
	\int_{z_\mathrm{s}-\epsilon}^{z_\mathrm{s}+\epsilon}\left[\rho(z)\frac{\mathrm{d}}{\mathrm{d}z}\left(\frac{1}{\rho(z)}\frac{\mathrm{d}\Psi}{\mathrm{d} z}\right)+\left(k^{2}-k_{r}^{2}\right)\Psi\right] \mathrm{d} z=-\frac{1}{2\pi}
\end{equation}
Since $\epsilon\rightarrow 0$, the above equation translates to:
\begin{equation}
	\label{eq.13}
	\left.\frac{\mathrm{d}\Psi}{\mathrm{d} z}\right|_{z_\mathrm{s}-\epsilon}^{z_\mathrm{s}+\epsilon}=-\frac{1}{2\pi} \Longleftrightarrow \left.\frac{\mathrm{d}\Psi}{\mathrm{d} z}\right|_{z_\mathrm{s}^+} - \left.\frac{\mathrm{d}\Psi}{\mathrm{d} z}\right|_{z_\mathrm{s}^-} =-\frac{1}{2\pi}
\end{equation}
This is the interface condition that the displacement potential at the depth of the sound source needs to satisfy. 

\subsection{Wavenumber integration}
Once the Green function of Eq.~\eqref{eq.4} or \eqref{eq.7} has been obtained, the corresponding displacement potential field can be obtained by evaluating the inverse Hankel/Fourier transform. However, in the actual numerical calculation, when using Eq.~\eqref{eq.3a} or \eqref{eq.6a} to compute the displacement potential field of a point source or a line source, only the finite interval $[k_{\min},k_{\max}]$ and $M$ discrete points can be used for numerical integration, where $k_{\min}$ and $k_{\max}$ are the lower and upper limits of numerical integration, respectively. Undersampling the spikes of the Green function with a limited number of discrete points would introduce large errors. In addition, acoustic waveguides have poles on or close to the real wavenumber axis. Fortunately, the aliasing problem can be eliminated by simply moving the integral contour out into the complex plane (see Chapter 4.5.5 of Ref.~\cite{Jensen2011}). For this purpose, a contour offset $\varepsilon$ can be introduced. If the points are chosen where the kernel functions are small and the contour offset satisfies $\varepsilon \ll (k_{\max}-k_{\min})$, then the contributions of the vertical parts become negligible compared to the integral along the horizontal part. Substituting $\bar{k}=k_r-\mathrm{i} \varepsilon$ into Eqs.~\eqref{eq.3a} and \eqref{eq.6a} can yield:

\begin{subequations}
	\label{eq.14}
	\begin{align}
		\label{eq.14a}
		&\psi(r,z)=\int_{0}^{\infty} \Psi\left(k_{r}-\mathrm{i} \varepsilon, z\right) J_{0}\left[(k_{r}-\mathrm{i} \varepsilon) r\right] \left(k_{r}-\mathrm{i} \varepsilon\right) \mathrm{d} k_{r}\\
		\label{eq.14b}
		&
		\begin{aligned}
			\psi(x, z)&=\int_{-\infty}^{\infty} \Psi\left(k_{x}-\mathrm{i} \varepsilon, z\right) \exp\left[{\mathrm{i} \left(k_{x}-\mathrm{i} \varepsilon\right) x}\right] \mathrm{d} k_{x}\\
			&=2\int_{0}^{\infty} \Psi\left(k_{x}-\mathrm{i} \varepsilon, z\right) \cos\left[{\left(k_{x}-\mathrm{i} \varepsilon\right) x}\right] \mathrm{d} k_{x}
		\end{aligned}
	\end{align}
\end{subequations}
For most practical purposes, the value of $\varepsilon$ can be taken to be \cite{Jensen2011}:
\begin{equation}
	\label{eq.15}
	\varepsilon=\frac{3\Delta k_r}{2 \pi \log _{10} \mathrm{e}},\quad \Delta k_r=\frac{k_{\max}-k_{\min}}{M-1}
\end{equation}
Therefore, the integrals in Eq.~\eqref{eq.14} become the rectangular integrals of the following form in the actual numerical calculations:

\begin{subequations}
	\label{eq.16}
	\begin{align}
		&\psi(r,z)=\Delta k_r \sum_{k_r=k_{\min}}^{k_{\max}}\Psi(k_r-\mathrm{i} \varepsilon,z) J_{0}\left[\left(k_{r}-\mathrm{i} \varepsilon\right) r\right] \left(k_{r}-\mathrm{i} \varepsilon\right)\\
		&\psi(x,z)=2\Delta k_x \sum_{k_x=k_{\min}}^{k_{\max}}\Psi(k_x-\mathrm{i} \varepsilon,z) \cos\left[\left(k_{x}-\mathrm{i} \varepsilon\right) x\right]
	\end{align}
\end{subequations}
$\Delta k_x$ has the same form as $\Delta k_r$ in Eq.~\eqref{eq.15}. The above numerical integration can be easily written in the form of matrix multiplication. In addition, we should choose the parameters of numerical integration with great care, such as $\Delta r$, $k_{\min}$, $k_{\max}$ and the farthest range $r_{\max}$ of the sound field of interest. Techniques related to evaluating integrals (e.g., quadrature schemes, fast field techniques) are universal; they are not the main innovations of this article and therefore are not described in detail.

When the displacement potential field is obtained by the above numerical integration, the sound pressure field can be obtained by the following formula \cite{Luowy2016}:

\begin{subequations}
	\label{eq.17}
	\begin{gather}
		p(r,z)=\rho(z) \omega^{2} \psi(r,z)\\
		p(x,z)=\rho(z) \omega^{2} \psi(x,z)
	\end{gather}
\end{subequations}
where $\omega = 2\pi f$. The transmission loss (TL) is defined as
\begin{equation}
	\label{eq.18}
	\mathrm{TL}=-20 \log _{10}\left|\frac{p}{p_{0}}\right|,\quad p_{0}=\begin{cases}
		\frac{\rho_\mathrm{s} \omega^{2}}{4 \pi}\exp({\mathrm{i} k_\mathrm{s}}), \quad \text{point source}\\
		\frac{\mathrm{i}\rho_\mathrm{s} \omega^{2}}{4} 	\mathcal{H}_{0}^{(1)}\left(k_\mathrm{s}\right), \quad \text{line source}
	\end{cases}	
\end{equation}
Here, $p_0$ is the acoustic pressure 1 m from the source, $\rho_\mathrm{s}$ and $k_\mathrm{s}$ are the density and wavenumber of the medium at the location of the source, respectively, and $\mathcal{H}_{0}^{(1)}(\cdot)$ denotes the first type of Hankel function.

From the above theory, we know that the Green function of Eq.~\eqref{eq.4} and \eqref{eq.7} is the key to the wavenumber integral model, and the main innovation of this paper is to propose a new algorithm based on spectral method to numerically solve the Green function of Eq.~\eqref{eq.4} and \eqref{eq.7}.

\section{Numerical Discretization}
\subsection{Chebyshev--Tau spectral method}
Next, we employ the Chebyshev--Tau spectral method to solve the depth-separated wave equation, Eq.~\eqref{eq.4}. The Chebyshev spectral method is a spectral method that uses Chebyshev polynomials as the basis functions \cite{Boyd2001}, so it is necessary to introduce these polynomials here.
\begin{equation}
	\begin{gathered}
		T_0(t)=1,\quad T_1(t)=t  \\
		T_{i+1}(t)=2tT_i(t)-T_{i-1}(t), \quad i\ge1
	\end{gathered}
\end{equation}
Chebyshev polynomials are a class of orthogonal polynomials whose orthogonality is defined as follows \cite{Jshen2011}:
\begin{equation}
	\int_{-1}^1 \frac{T_i(t)T_j(t)}{\sqrt{1-t^2}}\mathrm {d}t=\begin{cases}
		0,&i\neq j\\
		\pi,&i=j=0\\
		\frac{\pi}{2},&i=j\ge 1\\
	\end{cases}
\end{equation}
Since the Chebyshev polynomials $\{T_i(t)\}$, that is, the basis functions, are defined in $t\in[-1,1]$, the equation to be solved, Eq.~\eqref{eq.4}, must first be scaled to $t\in[-1,1]$ as follows:
\begin{equation}
	\label{eq.21}
	\mathcal{L}\Psi(t) = 0,\quad \mathcal{L}=
	\frac{4}{\vert \Delta h\vert^2}\rho(t)\frac{\mathrm{d}}{\mathrm{d}t}\left(\frac{1}{\rho(t)}\frac{\mathrm{d}}{\mathrm {d}t}\right) +\left[k^{2}(t)-k_r^2\right]
\end{equation}
where $\Delta h$ denotes the length of the domain and $\mathcal{L}$ represents the differential operator.

Next, the function to be determined, $\Psi(t)$, is transformed into the spectral space spanned by the basis functions $\{T_i(t)\}_{i=0}^\infty$. Furthermore, the expression for the spectral coefficients $\{\hat{\Psi}_i\}_{i=0}^\infty$ can also be obtained from the orthogonality of the Chebyshev polynomials \cite{Canuto2006}.
\begin{equation}
	\label{eq.22}
	\Psi(t)=\sum_{i=0}^{\infty}\hat{\Psi}_{i}T_i(t) \Longleftrightarrow 
	\hat{\Psi}_{i}=\frac{2}{\pi d_i}\int_{-1}^1 \frac{T_i(t)\Psi(t)}{\sqrt{1-t^2}}\mathrm {d}t,\quad d_i=\begin{cases}
		2,\quad i=0\\
		1,\quad i>0
	\end{cases}
\end{equation}
The integral on the right side of the above equation is usually calculated using the Gauss--Chebyshev--Lobatto numerical quadrature \cite{Boyd2001}.

Since it is impossible to expand to infinite terms in actual calculation, only the first $(N+1)$ terms can be retained \cite{Guoby1998}:
\begin{equation}
	\label{eq.23}
	\Psi(t)\approx \Psi_N(t)=\sum_{i=0}^{N}\hat{\Psi}_{i}T_i(t)
\end{equation}
$\Psi_N(t)$ is a function approximation, which becomes increasingly accurate as $N$ increases. The truncation of the infinite term expansion described above inevitably introduces errors, which means that Eq.~\eqref{eq.21} no longer strictly holds. Substituting $\Psi_N(t)$ into Eq.~\eqref{eq.21} yields a residual \cite{Gottlieb1977}, which we call $R_N(t)$.
\begin{equation}
	R_N(t)=\mathcal{L}\Psi_N(t)
\end{equation}
Some principle must be adopted to minimize $R_N$ so that the above spectral expansion can achieve higher accuracy. In a Tau-type spectral method, the basis functions are used as the weight functions, and then the inner products of the weight functions and the residuals are forced to equal 0 \cite{Lanczos1938}.
\begin{equation}
	\int_{-1}^1 \frac{\mathcal{L}\Psi_N(t)T_i(t)}{\sqrt{1-t^2}}\mathrm{d} t = 0,\quad i=0,1,\cdots,N
\end{equation}
This constraint on the residuals is the essence of the weighted residual method \cite{Jshen2011}. In mathematical monographs, the above equation is generally called the weak form or variational form of Eq.~\eqref{eq.21} \cite{Boyd2001}. Taking into account the orthogonality of the Chebyshev polynomials and Eq.~\eqref{eq.22}, the above equation becomes:
\begin{equation}
	\label{eq.26}
	\hat{\mathcal{L}}\hat{\Psi}_{i} = 0,\quad i=0,1,\cdots,N
\end{equation}
where $\hat{\mathcal{L}}$ represents the $\mathcal{L}$ operator on the spectral space.

The next most important thing is $\hat{\mathcal{L}}$, that is, the transformation of the $\mathcal{L}$ operator to the spectral space. The $\mathcal{L}$ operator contains a derivative term. According to the properties of the Chebyshev polynomials, the following is straightforward to prove:
\begin{equation}
	\label{eq.27}
	\hat{\Psi}'_i = \frac{2}{c_i}
	\sum_{\substack{j=i+1,\\ 
			j+i=\text{odd}
	}}^{N} j \hat{\Psi}_j, \quad c_0=2,c_{i>1}=1
	\Longleftrightarrow \bm{\hat{\Psi}}' = \mathbf{D}_N \bm{\hat{\Psi}}
\end{equation}
Thus, the derivative term is transformed into a differential matrix $\mathbf{D}_N$, which is related only to the truncation order $N$ and is completely unrelated to $\Psi$. The term is obtained from the relationship between Chebyshev polynomials and their derivatives \cite{Canuto2006}.

The $\mathcal{L}$ operator also contains a product term, and the spectral transformation of the product of the two functions satisfies the following relationship:
\begin{equation}
	\label{eq.28}
	\widehat{(v\Psi)}_i \approx 
	\frac{1}{2} \sum_{m+n=i}^{N} \hat{\Psi}_m\hat{v}_n +
	\frac{1}{2} \sum_{\vert m-n\vert=i}^{N} \hat{\Psi}_m\hat{v}_n  \Longleftrightarrow  \widehat{\bm{(v\Psi)}} \approx \mathbf{C}_v \bm{\hat{\Psi}}
\end{equation}
where $v=v(t)$ is any continuous function on $t\in[-1,1]$ used, for example. Similarly, the relationship between the spectral coefficients of the product of two functions and the spectral coefficients of the individual functions is represented by a matrix $\mathbf{C}_v$, which is related only to $v$ and not to $\Psi$ \cite{Orszag1972,Boyd2001}.

\subsection{Discretization}
According to the above analysis, Eq.~\eqref{eq.21} is discretized into the following matrix--vector form in the Chebyshev spectral space:
\begin{equation}
	\label{eq.29}
	\left(\frac{4}{\vert\Delta h\vert^2}\mathbf{C}_{\rho}\mathbf{D}_{N}\mathbf{C}_{1/\rho}\mathbf{D}_{N}+\mathbf{C}_{k^2}-k_r^2\mathbf{E}_N\right)\bm{\hat{\Psi}}= \mathbf{0}
\end{equation}
where $\mathbf{E}_N$ is the $(N+1)$-order identity matrix. The above equation is equivalent to Eq.~\eqref{eq.26}, where $\bm{\hat{\Psi}}$ is a column vector composed of $\{\hat{\Psi}_i\}_{i=0}^N$. Eq.~\eqref{eq.29} is a set of linear equations, and the boundary conditions are not imposed at this time.

For the waveguide in Fig.~\ref{Figure1}, the depth-separated wave equation must be established in all discontinuous layers. A single set of basis functions cannot span all layers since the Chebyshev polynomials do not satisfy the continuity conditions at the interfaces $\{h_l\}_{l=1}^\ell$. Thus, we apply the domain-decomposition strategy \cite{Min2005} to Eq.~\eqref{eq.4} and split the complete domain into $\ell$ subintervals:
\begin{equation}
	\label{eq.30}
	\Psi_l(z)=\Psi_l(t) \approx\sum_{i=0}^{N_l}\hat{\Psi}_{l,i}T_i(t),\quad
	t=\frac{2z}{h_{l-1}-h_{l}}+\frac{h_l+h_{l-1}}{h_l-h_{l-1}},\quad z\in[h_{l-1}, h_l]
\end{equation}
$N_l$ is the spectral truncation order in the $l$-th layer, and $\{\hat{\Psi}_{l,i}\}_{i=0}^{N_l}$ is the spectral coefficient in the $l$-th layer. Similar to Eq.~\eqref{eq.29}, the depth-separated wave equation in the $l$-th layer can be discretized into the following matrix--vector form:
\begin{equation}
	\label{eq.31}
	\mathbf{A}_l\bm{\hat{\Psi}}_l=\bm{0},\quad \mathbf{A}_l= 
	\frac{4}{(h_l-h_{l-1})^2}\mathbf{C}_{\rho_l}\mathbf{D}_{N_l}\mathbf{C}_{1/\rho_l}\mathbf{D}_{N_l}+\mathbf{C}_{k_l^2}-k_r^2\mathbf{E}_{N_l}
\end{equation}
where $\mathbf{A}_l$ is a square matrix of order $(N_l+1)$ and $\bm{\hat{\Psi}}_l$ is a column vector composed of $\{\hat{\Psi}_{l,i}\}_{i=0}^{N_l}$. Since the interface conditions are related to the adjacent layers, a total of $\ell$ Eqs.~\eqref{eq.31} of $l=1,\cdots,\ell$ should be solved simultaneously, which is expressed as follows:
\begin{equation}
	\label{eq.32}
	\left[\begin{array}{cccc}
		\mathbf{A}_1&\mathbf{0}&\mathbf{0}&\mathbf{0}\\
		\mathbf{0}&\mathbf{A}_2&\mathbf{0}&\mathbf{0}\\
		\mathbf{0}&\mathbf{0}&\ddots&\mathbf{0}\\
		\mathbf{0}&\mathbf{0}&\mathbf{0}&\mathbf{A}_\ell\\
	\end{array}\right]
	\left[\begin{array}{c}
		\bm{\hat{\Psi}}_1\\
		\bm{\hat{\Psi}}_2\\
		\vdots\\
		\bm{\hat{\Psi}}_\ell\\
	\end{array}
	\right]=\left[
	\begin{array}{c}
		\bm{0}\\
		\bm{0}\\
		\vdots\\
		\bm{0}\\
	\end{array}\right]
\end{equation}
Note that when $z_\mathrm{s}$ is not on the interface, we set up a virtual interface for it as described in Eq.~\eqref{eq.13}. Eq.~\eqref{eq.32} that modifies to be satisfied for the two layers above and below the virtual interface can also be organized into block diagonal form, and the total number of layers becomes $(\ell+1)$ at this time.
\begin{equation}
	\label{eq.33}
	\left[\begin{array}{ccccc}
		\mathbf{A}_1&\mathbf{0}&\mathbf{0}&\mathbf{0}&\mathbf{0}\\
		\mathbf{0}&\ddots&\mathbf{0}&\mathbf{0}&\mathbf{0}\\
		\mathbf{0}&\mathbf{0}&\mathbf{A}_\mathrm{s}&\mathbf{0}&\mathbf{0}\\
		\mathbf{0}&\mathbf{0}&\mathbf{0}&\ddots&\mathbf{0}\\
		\mathbf{0}&\mathbf{0}&\mathbf{0}&\mathbf{0}&\mathbf{A}_{\ell+1}\\
	\end{array}\right]
	\left[\begin{array}{c}
		\bm{\hat{\Psi}}_1\\
		\vdots\\
		\bm{\hat{\Psi}}_\mathrm{s}\\
		\vdots\\
		\bm{\hat{\Psi}}_{\ell+1}\\
	\end{array}
	\right]=\left[
	\begin{array}{c}
		\bm{0}\\
		\vdots\\
		\bm{0}\\
		\vdots\\
		\bm{0}\\
	\end{array}\right]
\end{equation}
where $\mathbf{A}_\mathrm{s}$ and $\bm{\hat{\Psi}}_\mathrm{s}$ represent Eq.~\eqref{eq.31} on the layer at the depth of the sound source.

The interface conditions in Eq.~\eqref{eq.8} and boundary conditions in Eqs.~\eqref{eq.9}--\eqref{eq.11} must also be expanded to the spectral space and explicitly added to Eq.~\eqref{eq.33}. In addition, on the virtual interface at the depth of the sound source, the intermittent condition, i.e., Eq.~\eqref{eq.13}, is also added to Eq.~\eqref{eq.33}. Upon considering the virtual interface, the seawater media comprise a total of $(\ell+1)$ layers, so there are $\ell$ interfaces leading to the $2\ell$ interface conditions. With the addition of the boundary conditions at the sea surface $(z=0)$ and seabed $(z=H)$, there are $2(\ell+1)$ conditions to apply. Next, we describe the imposition of boundary conditions and interface conditions in detail. For convenience of description, we define the following intermediate row vectors:
\[
\mathbf{s}_l=[s_0,s_1,s_2,\dots,s_{N_l}],\quad \mathbf{q}_l=[q_0,q_1,q_2,\dots,q_{N_l}]   
\]
where $s_i=T_i(-1)=(-1)^i$, $q_i=T_i(+1)=1$. Thus, the interface conditions and boundary conditions of Eqs.~\eqref{eq.8}--\eqref{eq.11} and \eqref{eq.13} can be transformed in Chebyshev spectral space as:

\begin{subequations}
	\begin{gather}
		\rho(h_l^+)\mathbf{s}_l\bm{\hat{\Psi}}_{l+1}-\rho(h_l^-)\mathbf{q}_l\bm{\hat{\Psi}}_l=0\\
		\mathbf{s}_l\mathbf{D}_{N_{l+1}+1}\bm{\hat{\Psi}}_{l+1}/(h_l-h_{l+1})-\mathbf{q}_l\mathbf{D}_{N_{l}+1}\bm{\hat{\Psi}}_l/(h_{l-1}-h_l)=0 \\
		\mathbf{q}_1\bm{\hat{\Psi}}_1=0  \\
		\mathbf{s}_{\ell+1}\bm{\hat{\Psi}}_{\ell+1}=0  \\
		\mathbf{s}_{\ell+1}\mathbf{D}_{N_{\ell+1}+1}\bm{\hat{\Psi}}_{\ell+1}=0  \\
		\left[2\rho_\infty \mathbf{s}_{\ell+1}\mathbf{D}_{N_{\ell+1}+1}/(h_l-h_{l+1}) + \rho(H)\sqrt{k_r^2-k_\infty^2} \mathbf{s}_{\ell+1}\right]\bm{\hat{\Psi}}_{\ell+1} = 0\\		
		\mathbf{s}_{\mathrm{s}+1}\mathbf{D}_{N_{\mathrm{s}+1}+1}\bm{\hat{\Psi}}_{\mathrm{s}+1}/(h_\mathrm{s}-z_\mathrm{s})-\mathbf{q}_{\mathrm{s}}\mathbf{D}_{N_\mathrm{s}+1}\bm{\hat{\Psi}}_\mathrm{s}/(z_\mathrm{s}-h_{\mathrm{s}+1})=-\frac{1}{4\pi}
	\end{gather}
\end{subequations}
where $h_\mathrm{s}$ and $h_{\mathrm{s}+1}$ represent the depth of the interfaces above and below the sound source, respectively.

How do these $2(\ell+1)$ conditions apply to Eq.~\eqref{eq.33}? A natural idea is to replace the last two rows of the $\mathbf{A}_1$ to $\mathbf{A}_{\ell+1}$ block matrix with the boundary or interface conditions that the corresponding layers need to satisfy. Doing so reduces the original spectral accuracy of each layer from order $N_l$ to order $(N_l-2)$, but this problem can be compensated by increasing the value of $N_l$. The spectral coefficients $\{\bm{\hat{\Psi}}_l\}_{l=1}^{\ell+1}$ of each layer of the Green function can be obtained by solving Eq.~\eqref{eq.33} after adding boundary constraints. The numerical solution of $\Psi(k_r,z)$ can be determined by performing the inverse Chebyshev transform (Eq.~\eqref{eq.23}) of $\{\bm{\hat{\Psi}}_l\}_{l=1}^{\ell+1}$ sequentially and then stacking into a single column vector.

\section{Numerical Simulation}
We present a program named WISpec (Wavenumber Integration based on the Spectral method) developed based on the above algorithm and verify the accuracy of the algorithm through several numerical experiments.

\subsection{Analytical example: ideal fluid waveguide}
The ideal fluid waveguide is a very simple example with an analytical solution. It consists of a layer of homogeneous seawater and upper/lower boundaries; the sea surface is usually perfectly free, and the bottom can be perfectly free or rigid. The ideal fluid waveguide of the perfectly free seabed has an analytical solution of the following form \cite{Jensen2011}:
\begin{equation}
	\label{eq.35}
	\begin{gathered}
		p(r,z)=\frac{2\pi \mathrm{i}}{H}\sum_{m=1}^{\infty}\sin{(k_{z,m}z_\mathrm{s})}\sin{(k_{z,m}z)}\mathcal{H}_0^{(1)}(k_{r,m}r)\\
		k_{z,m}=\frac{m\pi}{H},\quad k_{r,m}=\sqrt{k^2-k_{z,m}^2},\quad m=1,2,3,\cdots
	\end{gathered}
\end{equation}
The analytical solution of the sound field of the perfectly rigid seabed is the same as that of Eq.~\eqref{eq.35}, except that the vertical wavenumber becomes:
\begin{equation}
	k_{z,m}=\left(m-\frac{1}{2}\right)\frac{\pi}{H},\quad m=1,2,3,\cdots
\end{equation}

In this example, the frequency of the sound source is $f=20$ Hz, we take the sea depth $H=100$ m, $z_\mathrm{s}=36$ m, the density $\rho=1$ g/cm$^3$, the speed of sound $c=1500$ m/s, and the maximum horizontal range $r_{\max}=3000$ m. The number of discrete points in the wavenumber domain is taken as $M=2048$, the integral interval is $[0,2k_0]$ ($k_0$ is the wavenumber in water), and the spectral truncation order is $N=10$.

\begin{figure}[htbp]
	\centering
	\subfigure[]{\label{Figure2a}\includegraphics[width=8cm]{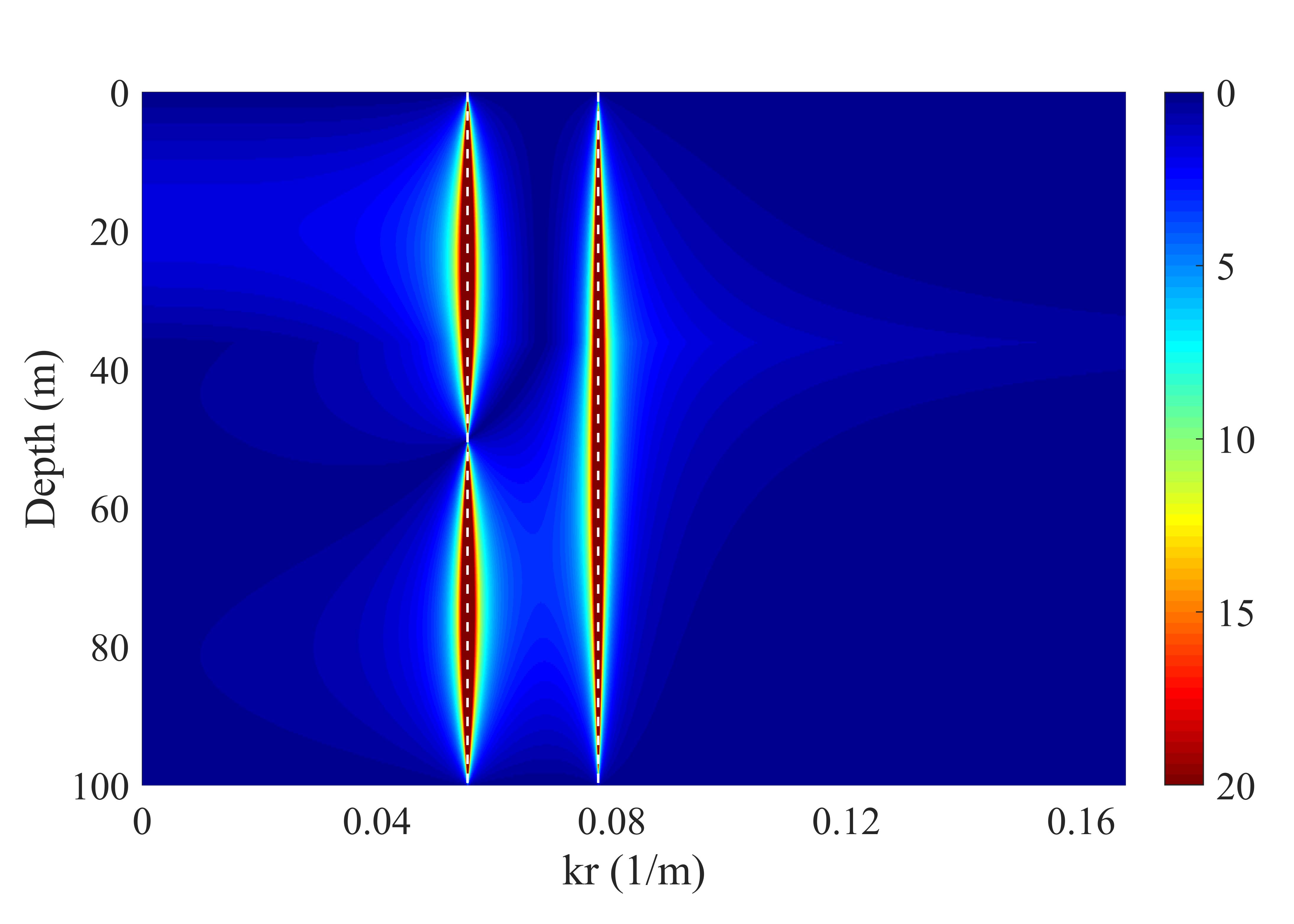}}
	\subfigure[]{\label{Figure2b}\includegraphics[width=8cm]{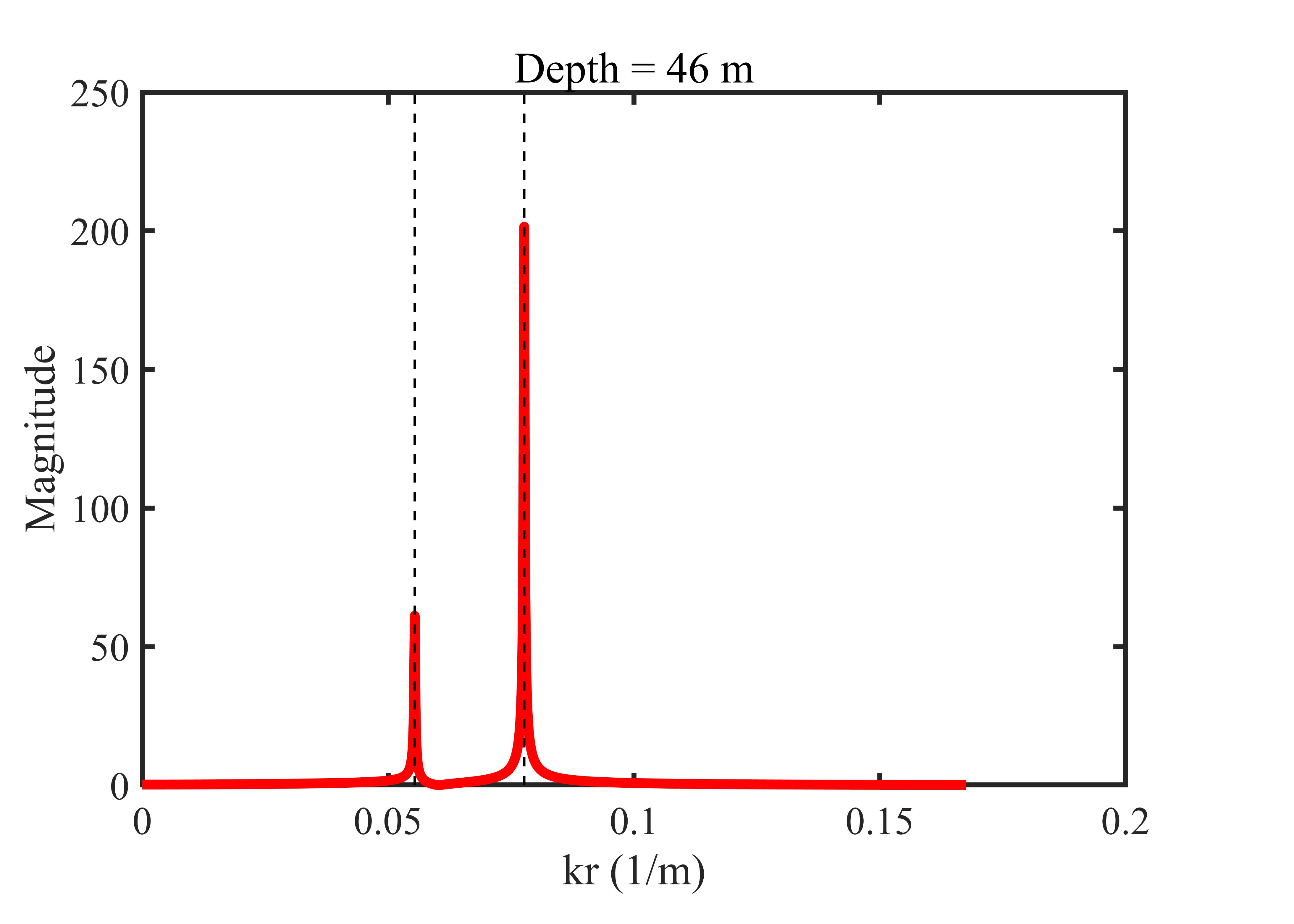}}\\
	\subfigure[]{\label{Figure2c}\includegraphics[width=8cm]{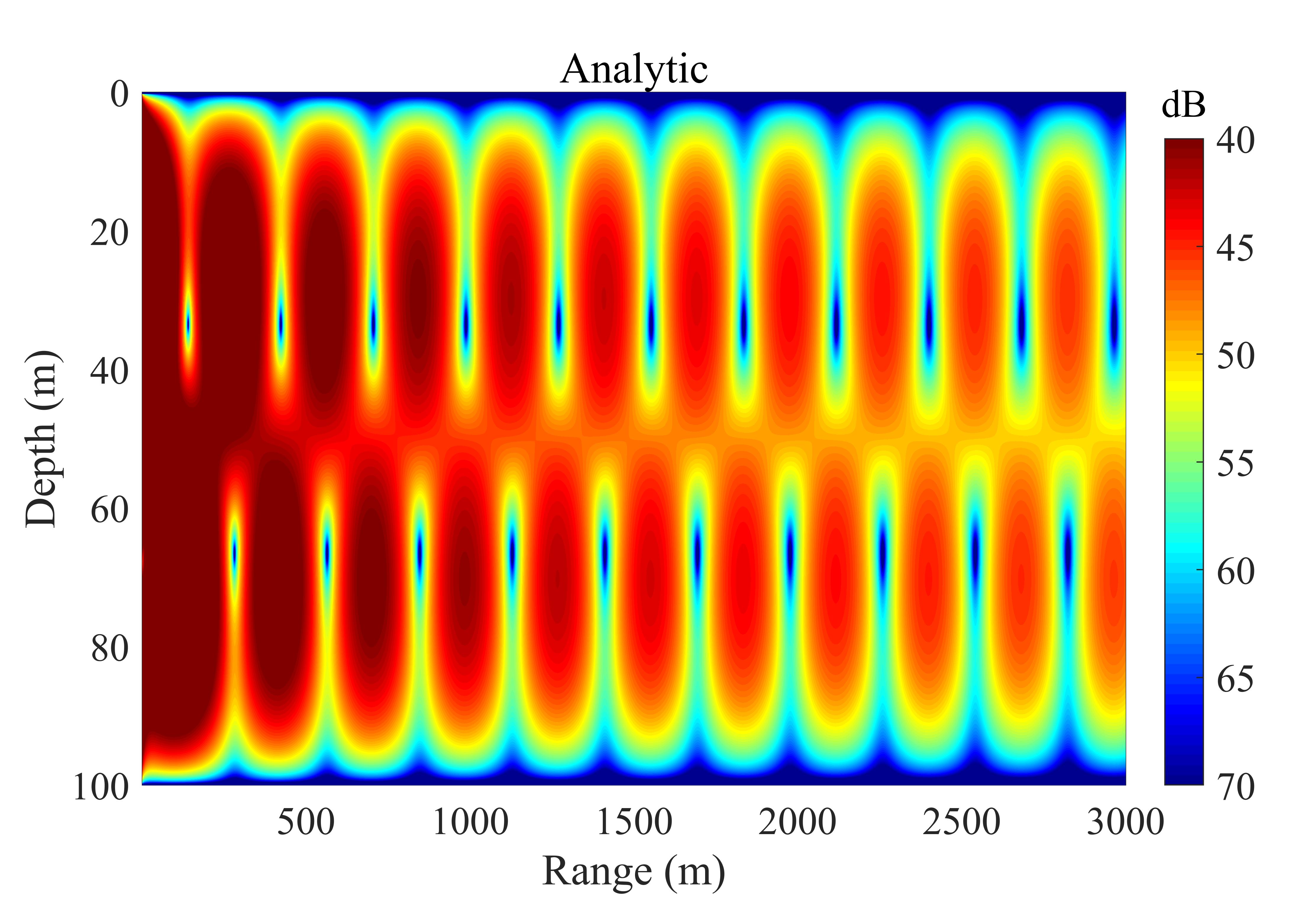}}
	\subfigure[]{\includegraphics[width=8cm]{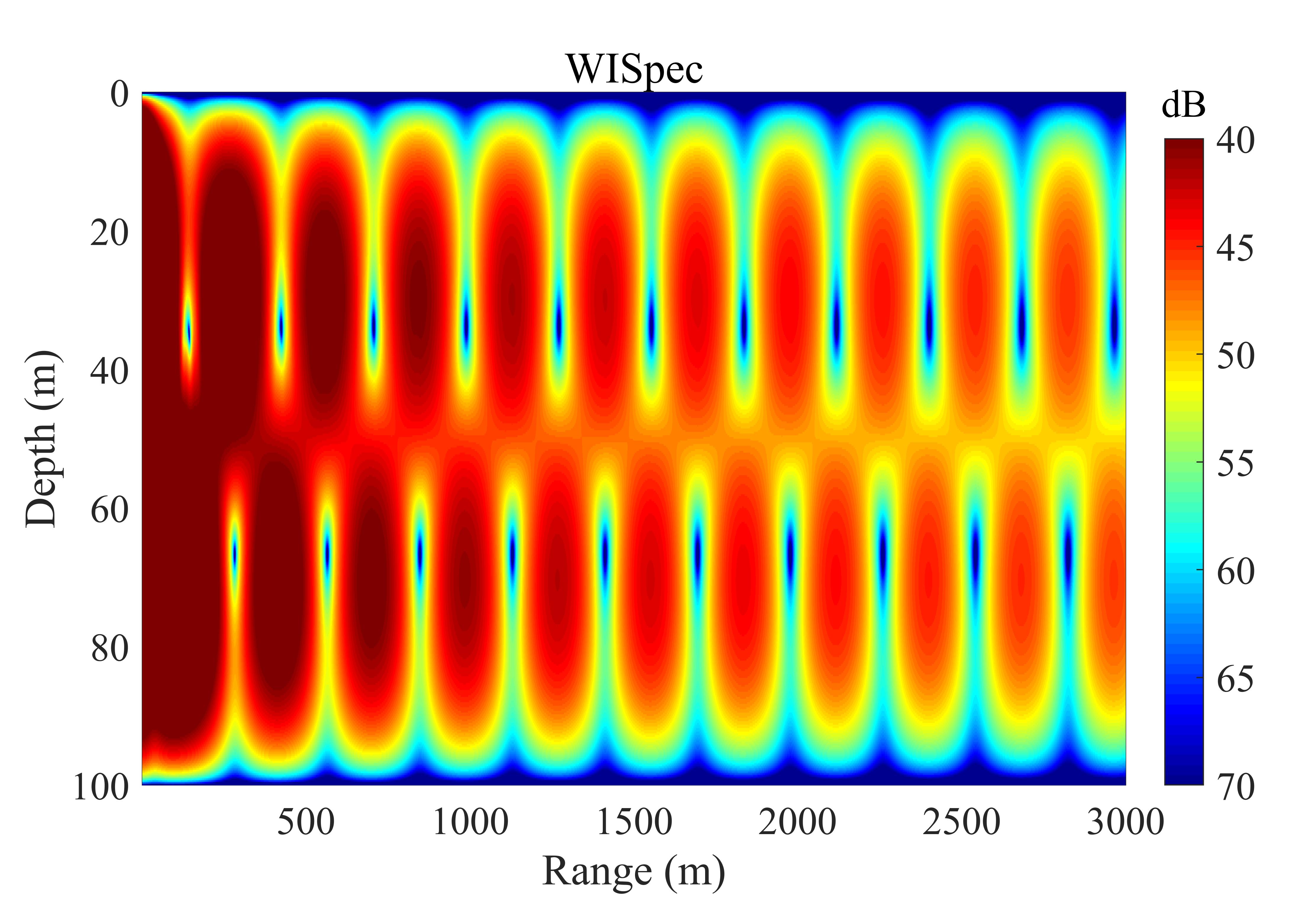}}\\
	\subfigure[]{\label{Figure2e}\includegraphics[width=8cm]{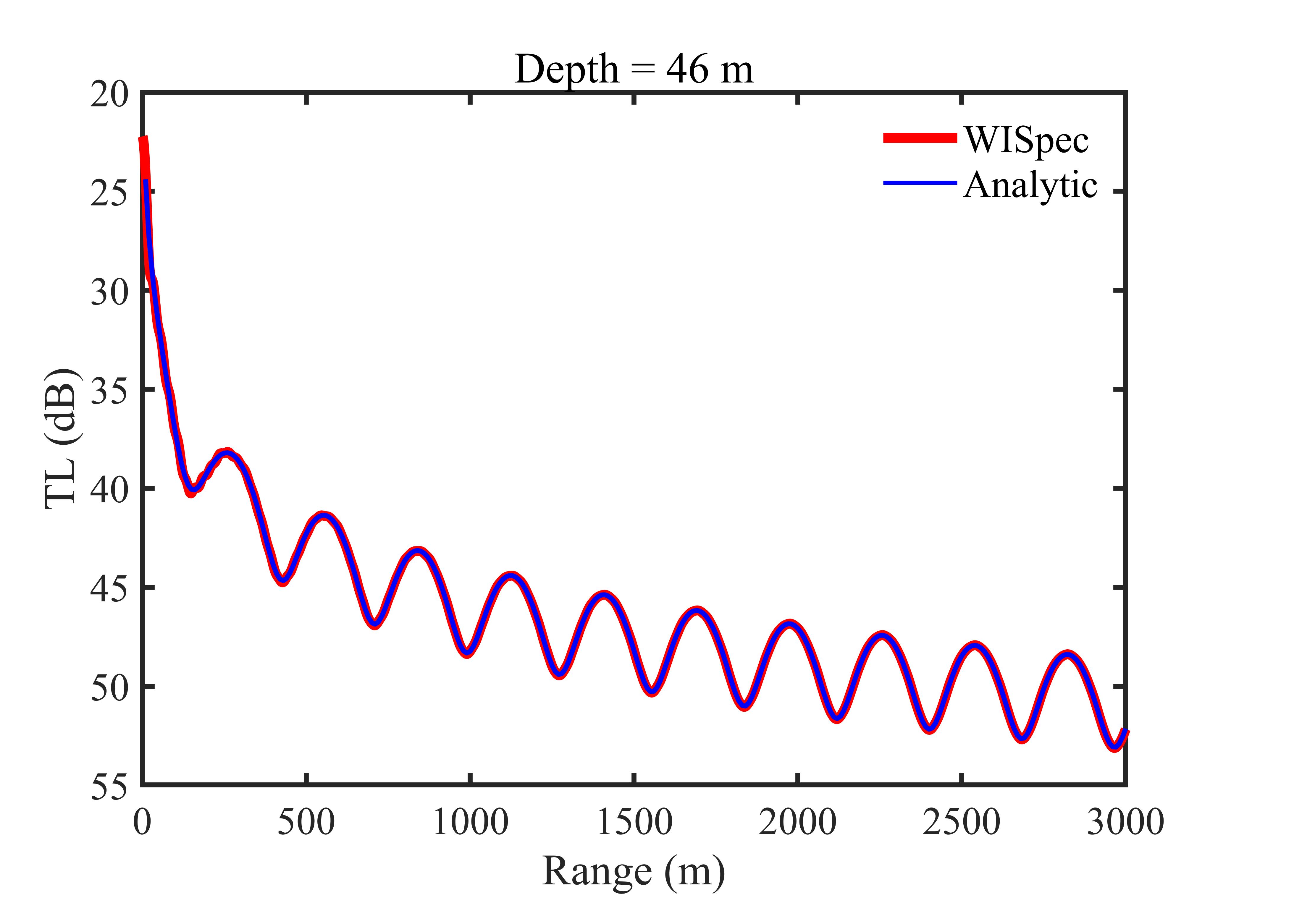}}
	\subfigure[]{\label{Figure2f}\includegraphics[width=8cm]{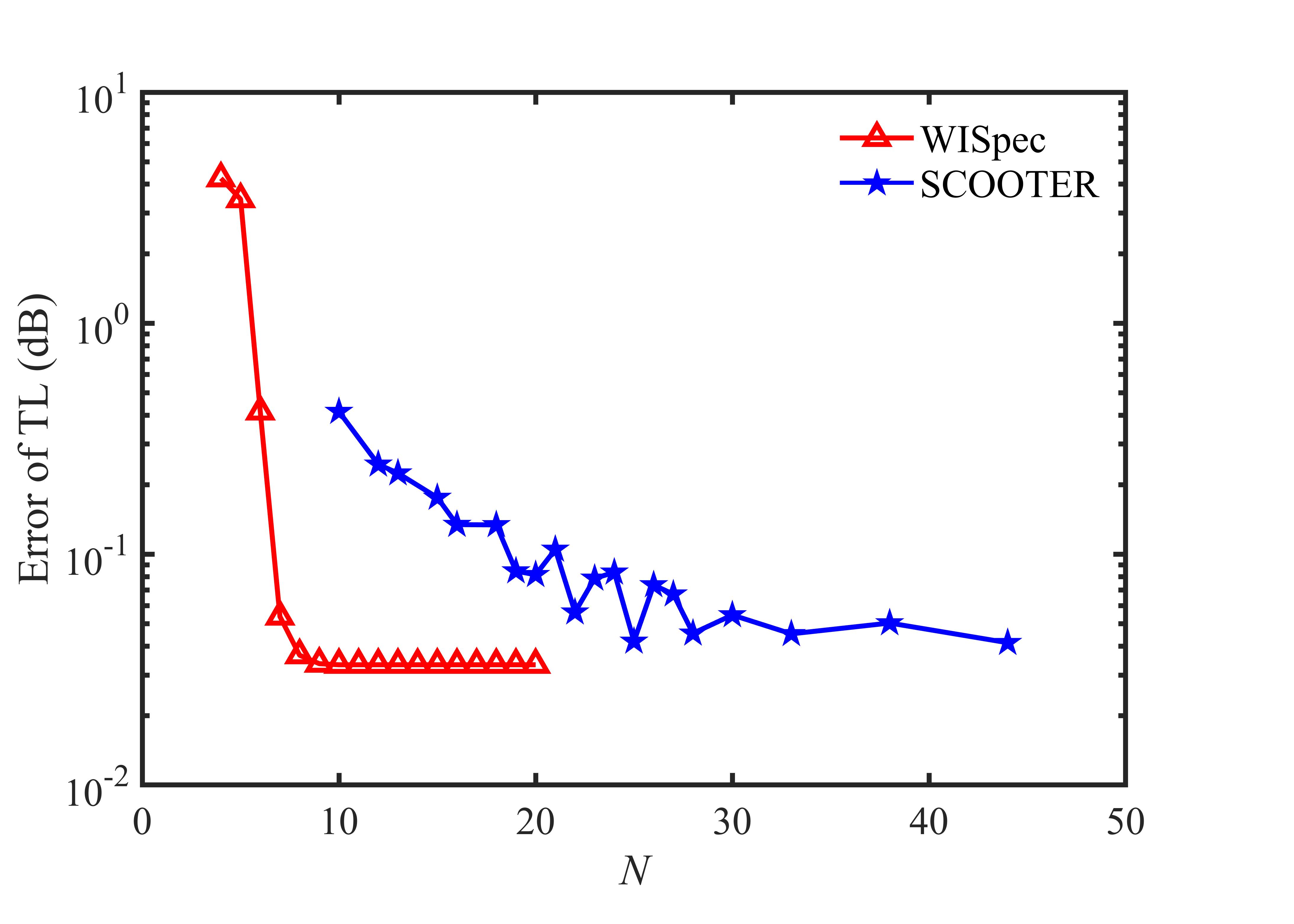}}
	\caption{Wavenumber spectrum of the ideal fluid waveguide with a perfectly free bottom calculated by WISpec (a); wavenumber spectrum at a depth of $z=46$ m (b); sound fields calculated via the analytical solution (c) and WISpec (d); TLs along the $r$-direction at a depth of $z=46$ m (e); and the variation trends with $N$ (f).}
	\label{Figure2}
\end{figure}

Fig.~\ref{Figure2a} and \ref{Figure2b} show the wavenumber spectrum of the ideal fluid waveguide with a perfectly free seabed calculated by WISpec. Two peaks (propagating modes) appear at the positions of $k=0.077678$ m$^{-1}$ and $0.055414$ m$^{-1}$. The presence of peaks in the wavenumber spectrum indicates that modes are excited by the source. The analytical solution of the discrete spectrum of the ideal fluid waveguide is listed in Table~\ref{tab1}. The wavenumbers at the peaks are very close to the analytical solution. At this time, the number of discrete points in the wavenumber domain is only 2048. The more discrete points there are in the wavenumber domain, the closer the wavenumbers at the peaks are to the analytical solution. Fig.~\ref{Figure2c} to \ref{Figure2e} show the pseudocolor and line graphs of the sound field calculated by the analytical solution and WISpec, respectively. The sound fields calculated by WISpec are very consistent with the analytical solution. Fig.~\ref{Figure2f} displays the variation trend of the error in the TL fields calculated by WISpec and SCOOTER with $N$. SCOOTER is a finite element code for simulating range-independent underwater acoustic fields \cite{SCOOTER}. In WISpec, $N$ is the spectral truncation order; in SCOOTER, $N$ is the number of piecewise linear elements, the same below. Here, to show the accuracy of the spectral method in solving the depth-separated wave equation, WISpec and SCOOTER use the same wavenumber integration parameters. Therefore, the error in the numerical sound fields almost completely comes from the error in the Green function. The numerical sound fields are given in the form of a discrete grid with 3000 discrete points horizontally from 1 to 3000 m and 401 discrete points vertically from 0 to 100 m. The error in the numerical sound field is calculated by:
\begin{equation}
	\text{TL}_{\text{error}}=\frac{\sum_{i=1}^{nz}\sum_{j=1}^{nr}\left|\text{TL}_{i,j}-\overline{\text{TL}}_{i,j}\right|}{nz\times nr}
\end{equation}
where $nz$ and $nr$ are the number of discrete points in the vertical and horizontal directions, respectively, and $\overline{\text{TL}}_{i,j}$ represents the analytical solution for the TL at $(z_i,r_j)$. Fig.~\ref{Figure2f} clearly illustrates that as $N$ increases, the error in the sound field rapidly converges to a very low level and remains stable, which also proves that the spectral method indeed maintains the advantage of exponential convergence in solving depth-separated wave equations. In Fig.~\ref{Figure2f}, the error associated with SCOOTER oscillates and decreases linearly, while WISpec quickly converges, thus further illustrating the efficiency of the spectral method.

\begin{table}[htbp]
	\centering
	\caption{Discrete modes of ideal fluid waveguides (unit: m$^{-1}$).}
		\begin{tabular}{ccccc}
			\hline
			\multirow{2}{*}{Mode order} &\multicolumn{2}{c}{Free seabed} &\multicolumn{2}{c}{Rigid seabed}   \\
			\cline{2-5}
			& Analytic & WISpec &Analytic & WISpec\\
			\hline
			1  &0.077 662 &0.077678 &0.082 290 &0.082343\\
			2  &0.055 412 &0.055414 &0.069 266 &0.069247\\
			3  &--        &--       &0.029 153 &0.029139\\
			\hline
		\end{tabular}
	\label{tab1}
\end{table}

\begin{figure}[htbp]
	\centering
	\subfigure[]{\label{Figure3a}\includegraphics[width=8cm]{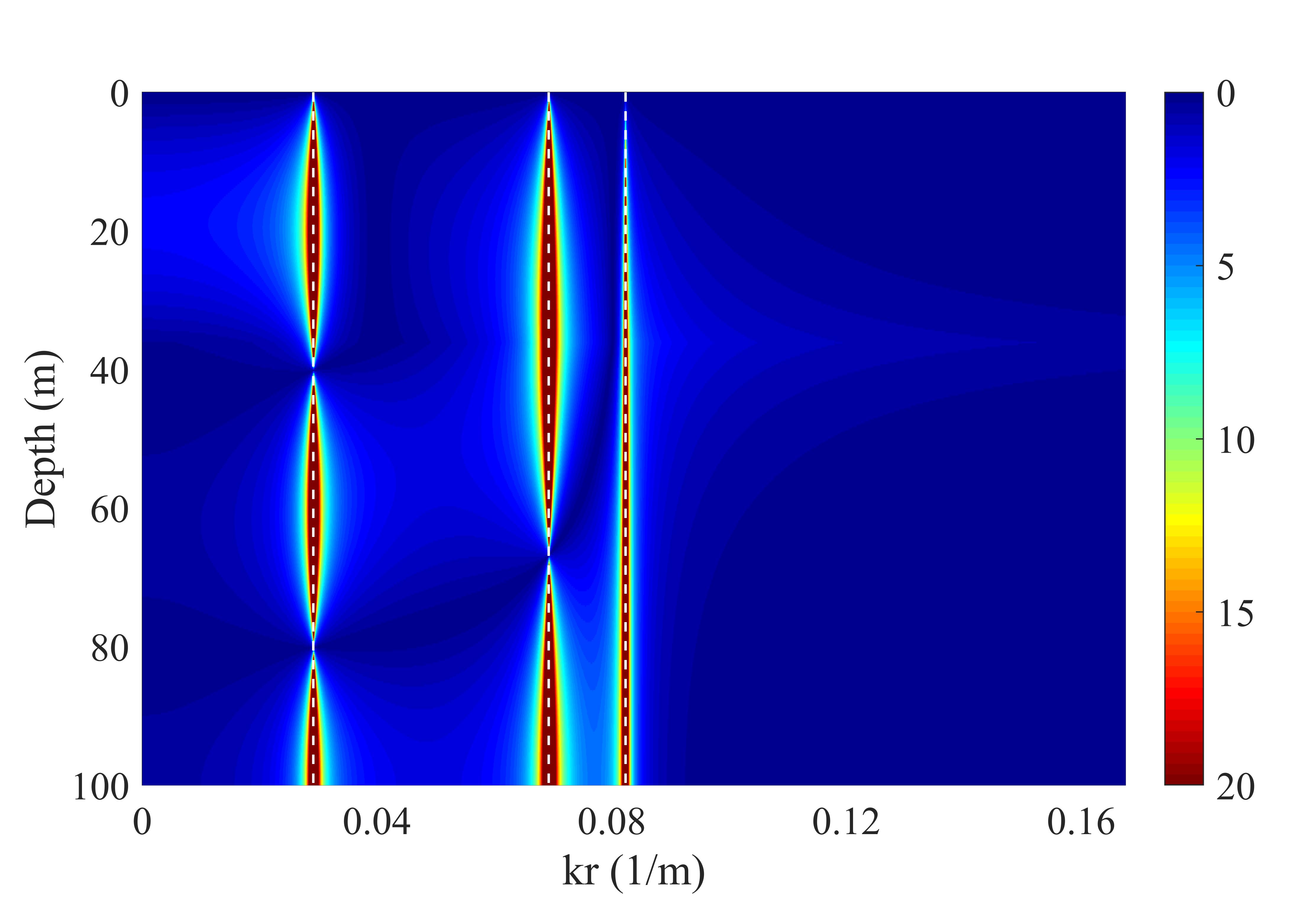}}
	\subfigure[]{\label{Figure3b}\includegraphics[width=8cm]{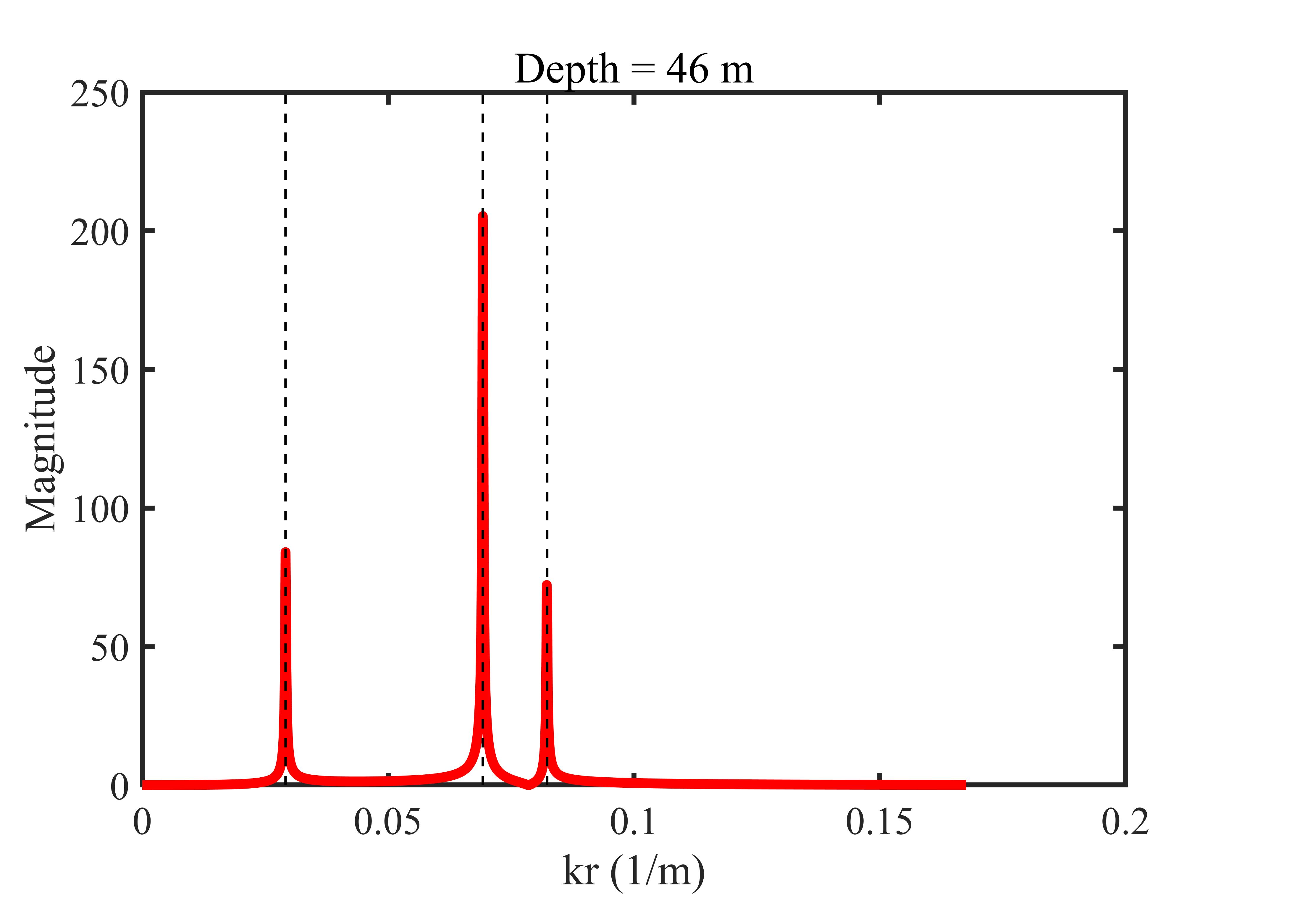}}\\
	\subfigure[]{\label{Figure3c}\includegraphics[width=8cm]{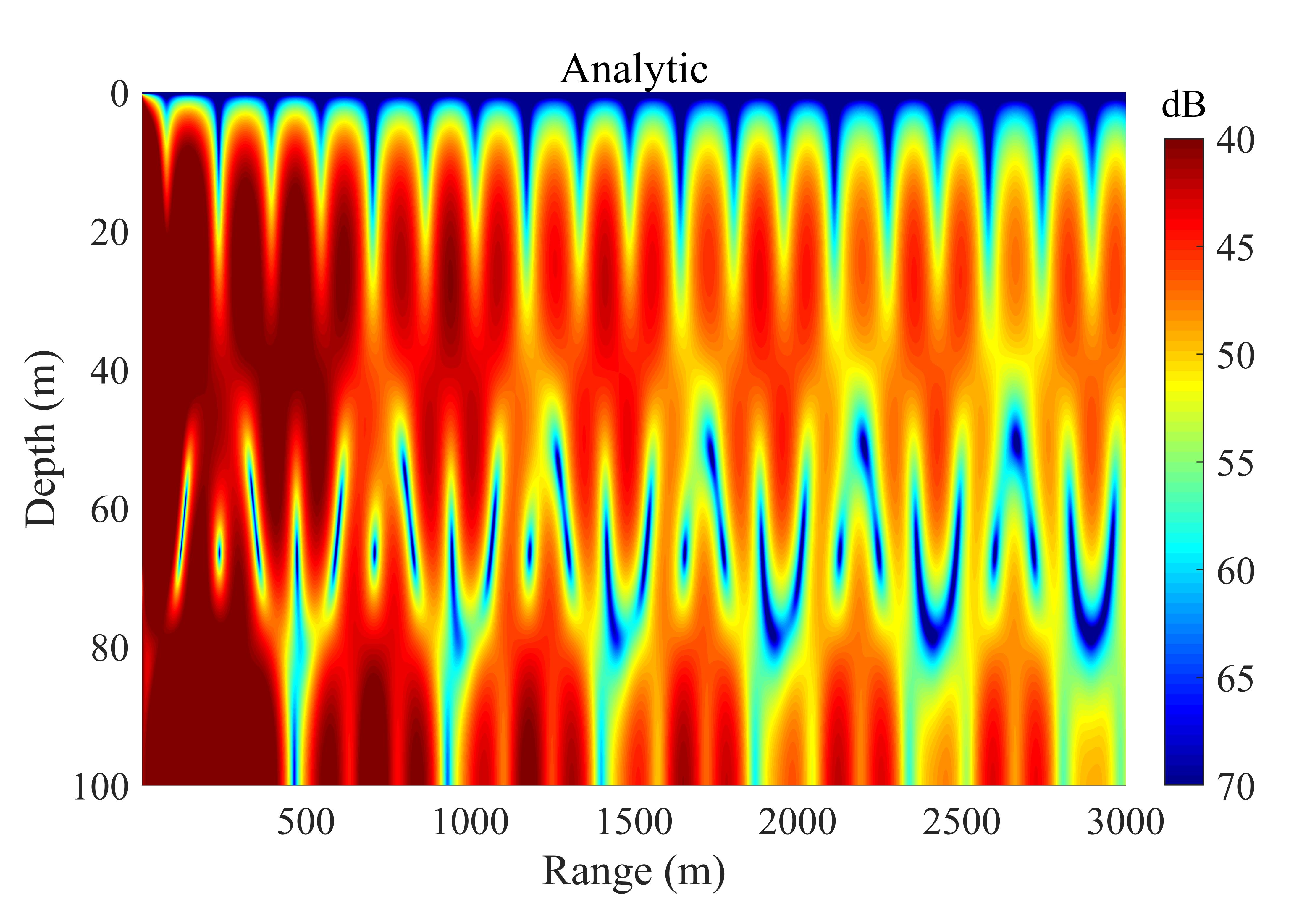}}
	\subfigure[]{\includegraphics[width=8cm]{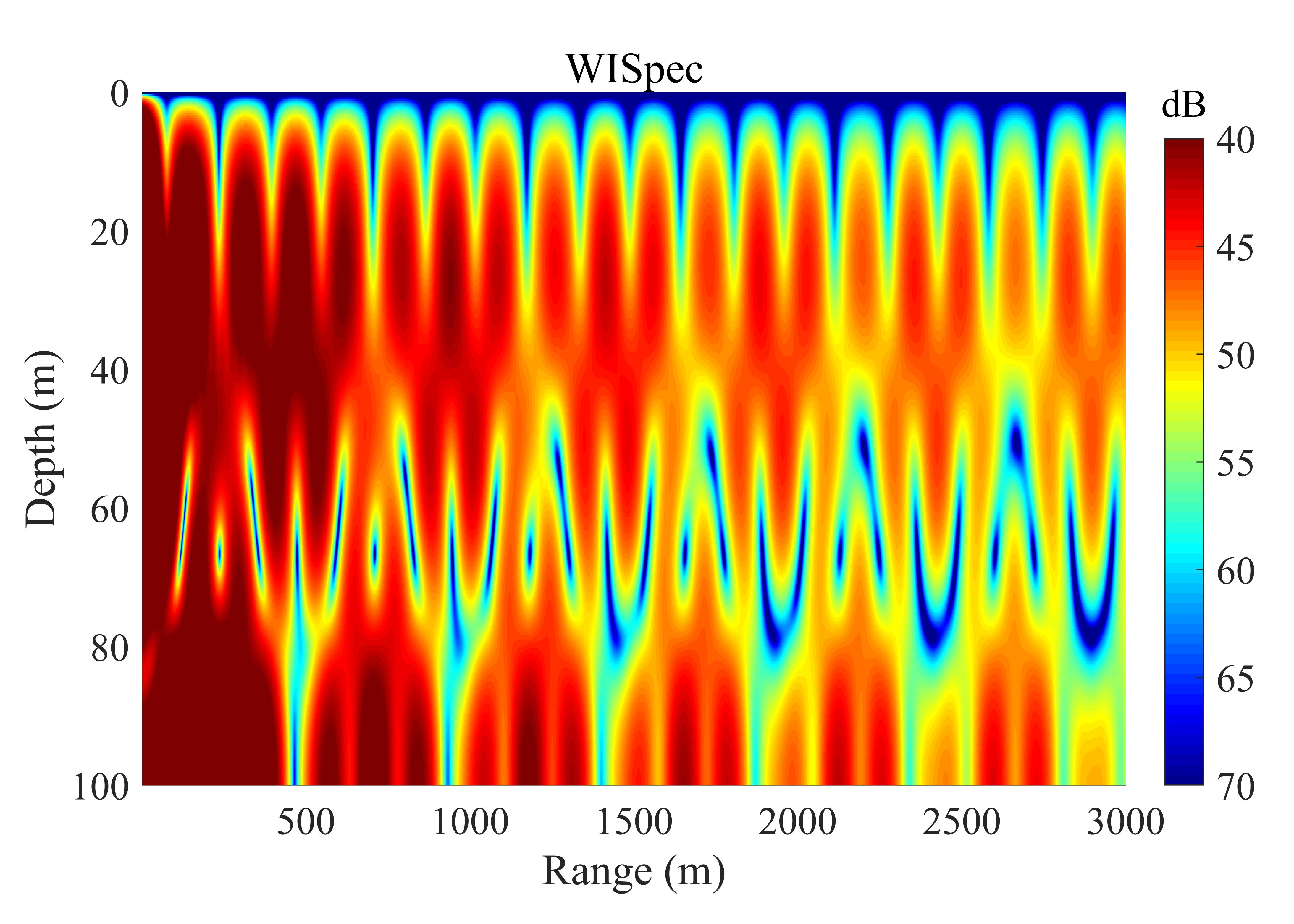}}\\
	\subfigure[]{\includegraphics[width=8cm]{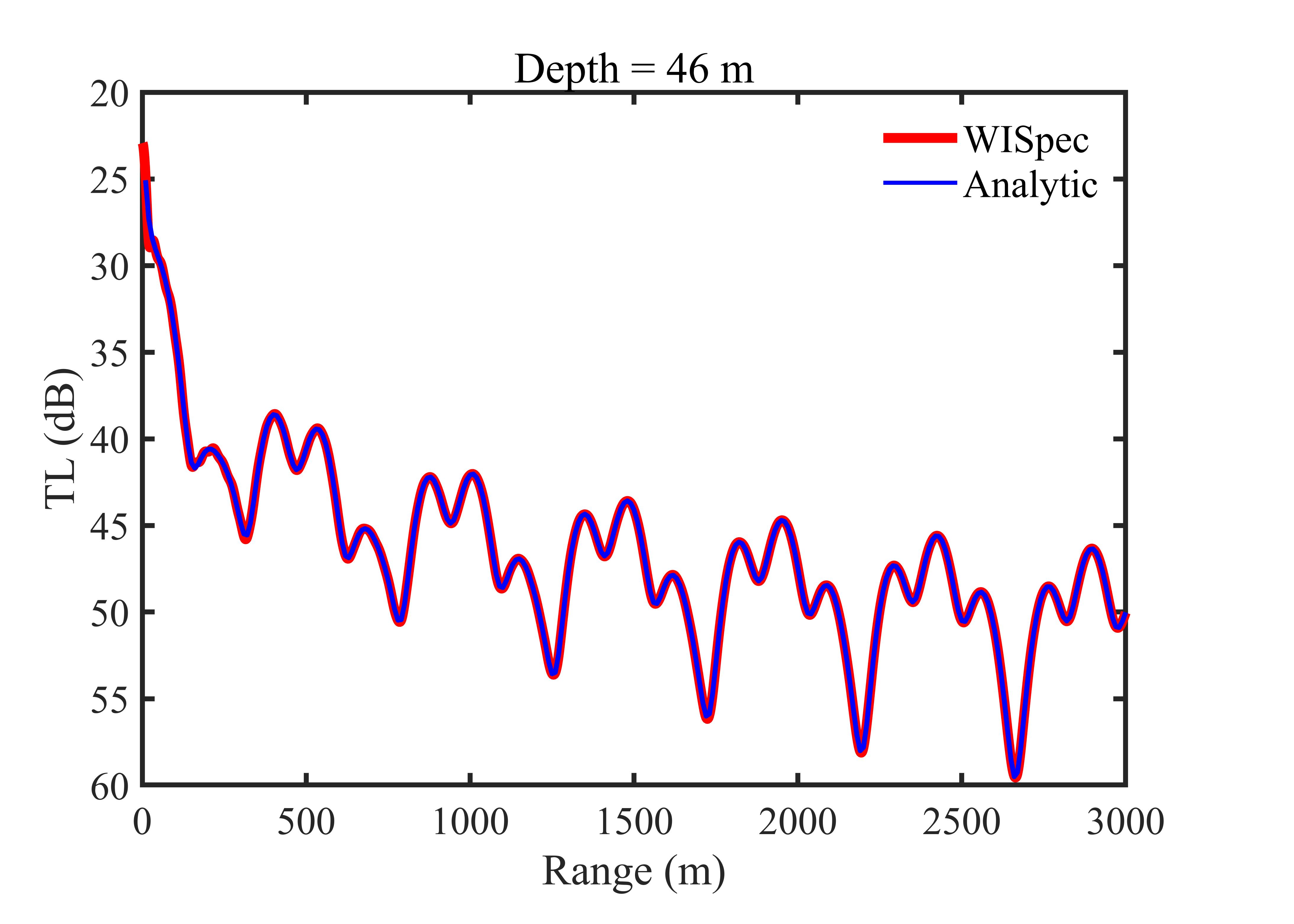}}
	\subfigure[]{\label{Figure3f}\includegraphics[width=8cm]{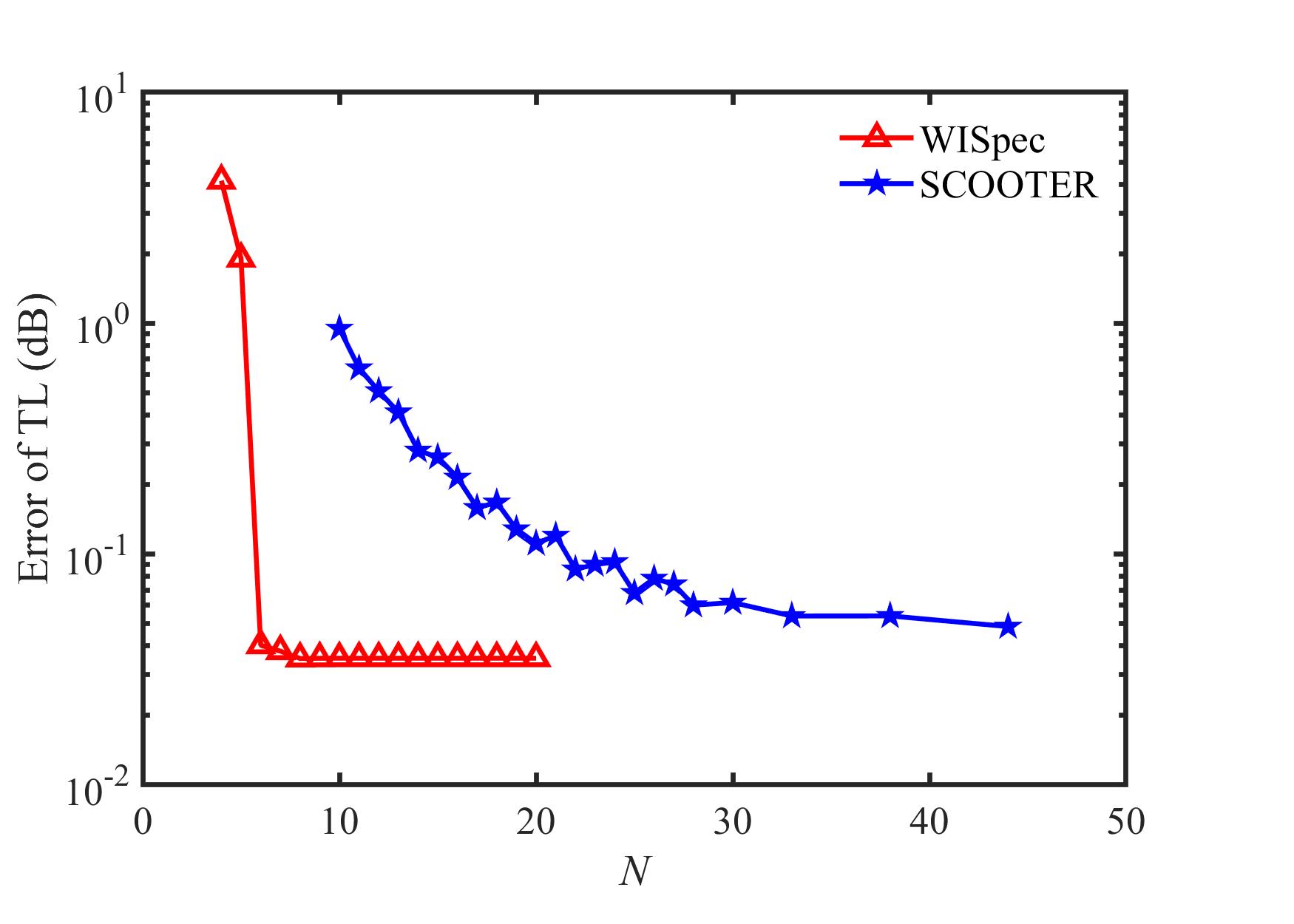}}
	\caption{Wavenumber spectrum of the ideal fluid waveguide with a perfectly rigid bottom calculated by WISpec (a); wavenumber spectrum at a depth of $z=46$ m (b); sound fields calculated via the analytical solution (c) and WISpec (d); TLs along the $r$-direction at a depth of $z=46$ m (e); and the variation trend of the error in the TL field with $N$ (f).}
	\label{Figure3}
\end{figure}

Similarly, Fig.~\ref{Figure3a} and \ref{Figure3b} show the wavenumber spectrum of the ideal fluid waveguide with the perfectly rigid seabed calculated by WISpec; three peaks appear at the positions of $k=0.082343$ m$^{-1}$, $0.069247$ m$^{-1}$ and $0.029139$ m$^{-1}$. This matches the analytical solution in Table~\ref{tab1} very well, and the sound fields shown in Fig.~\ref{Figure3c} to \ref{Figure3f} lead to the same conclusion as in Fig.~\ref{Figure2}, namely, that WISpec can calculate the sound field very accurately.

\subsection{Analytical example: pseudolinear-speed waveguide}
\begin{figure}[htbp]
	\centering
	\subfigure[]{\label{Figure4a}\includegraphics[width=8cm]{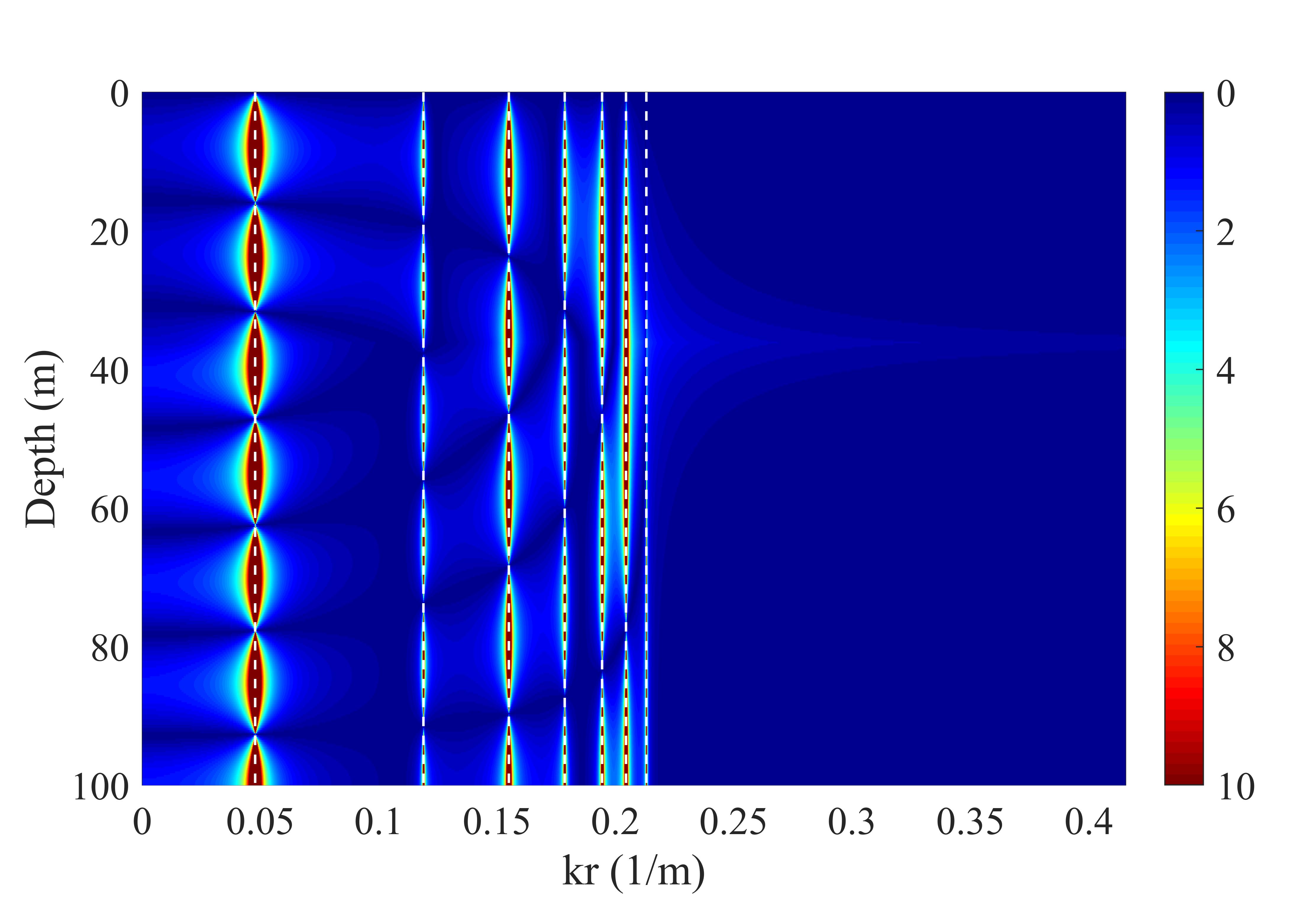}}
	\subfigure[]{\label{Figure4b}\includegraphics[width=8cm]{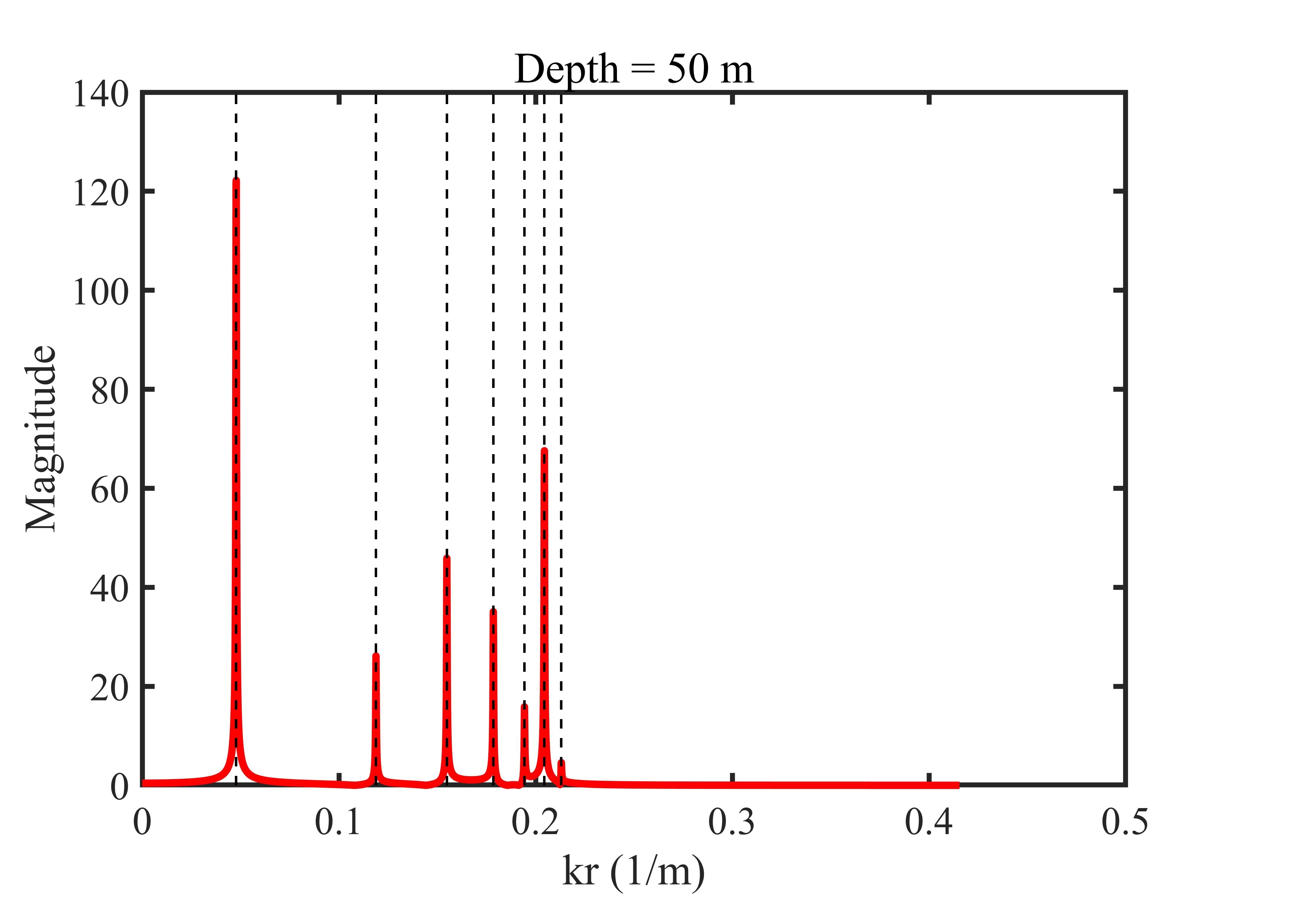}}\\
	\subfigure[]{\label{Figure4c}\includegraphics[width=8cm]{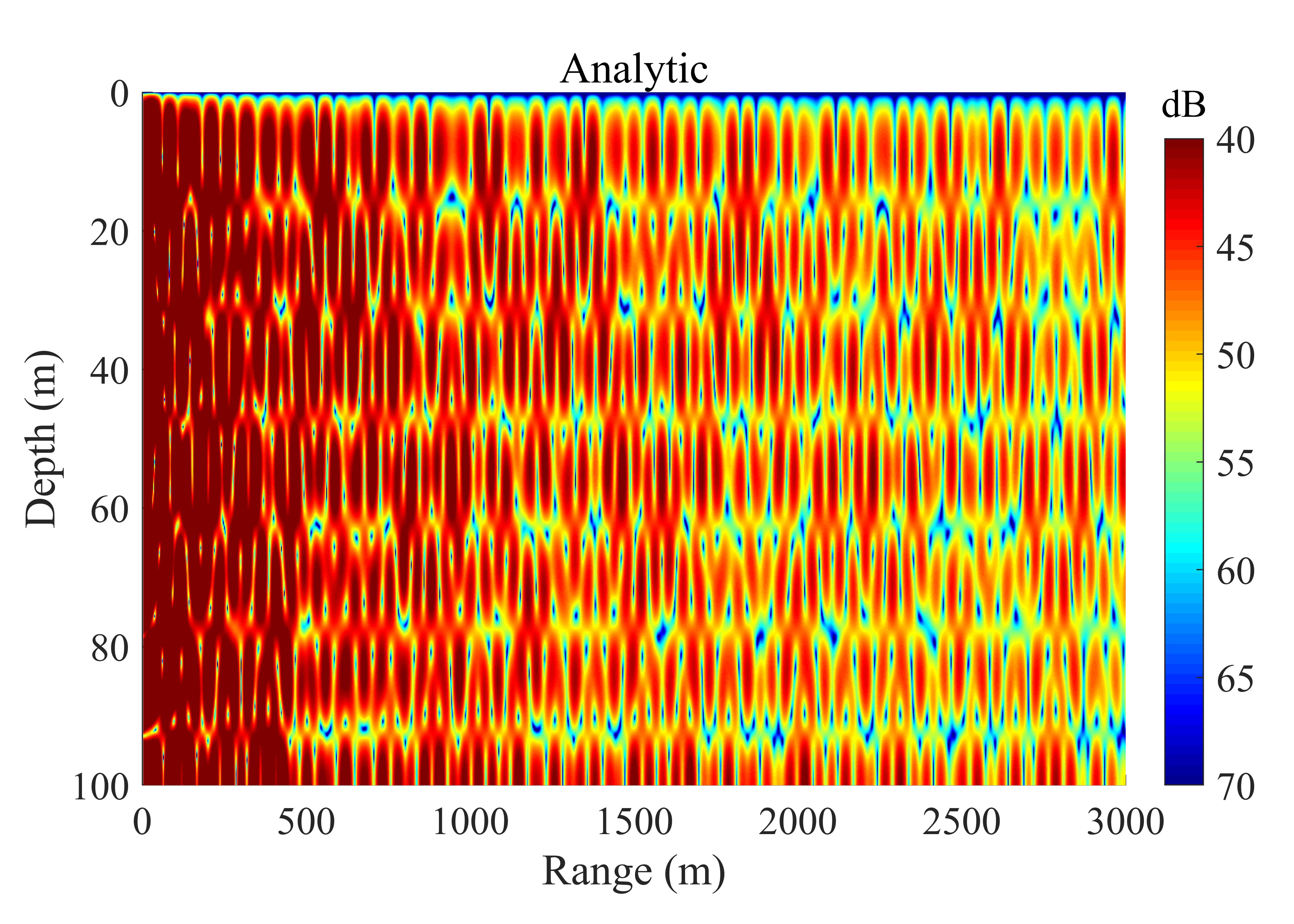}}
	\subfigure[]{\includegraphics[width=8cm]{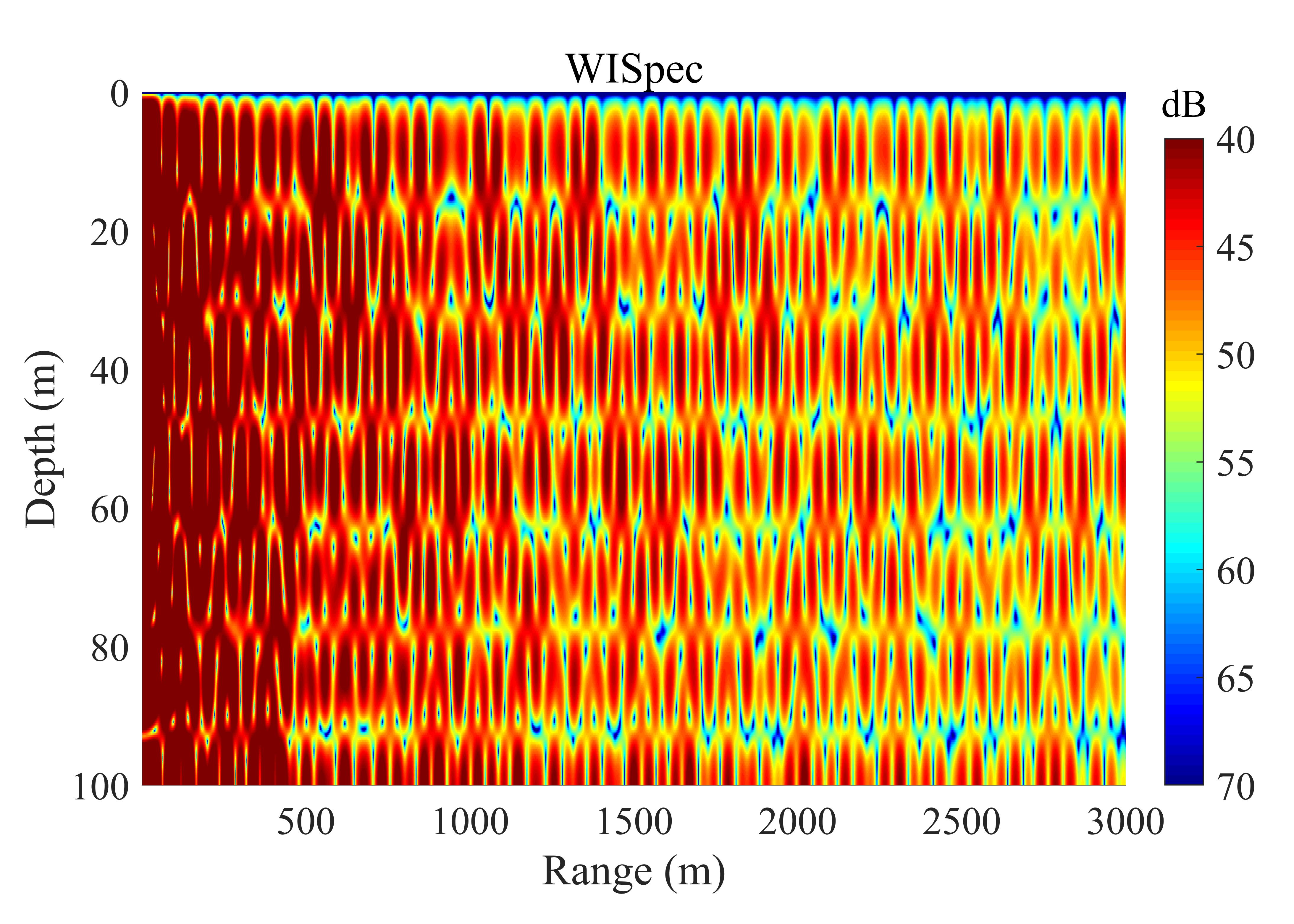}}\\
	\subfigure[]{\label{Figure4e}\includegraphics[width=8cm]{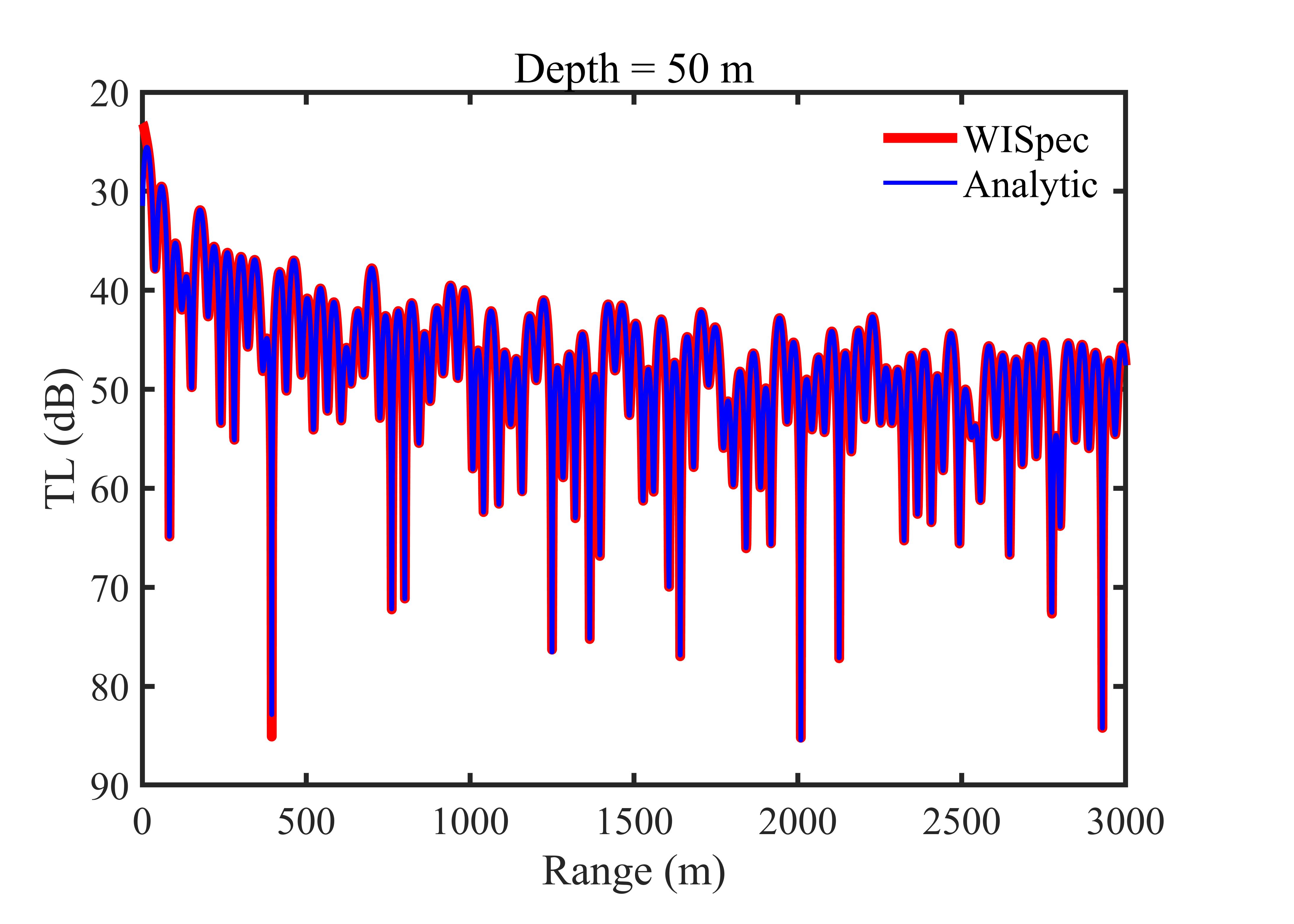}}
	\subfigure[]{\label{Figure4f}\includegraphics[width=8cm]{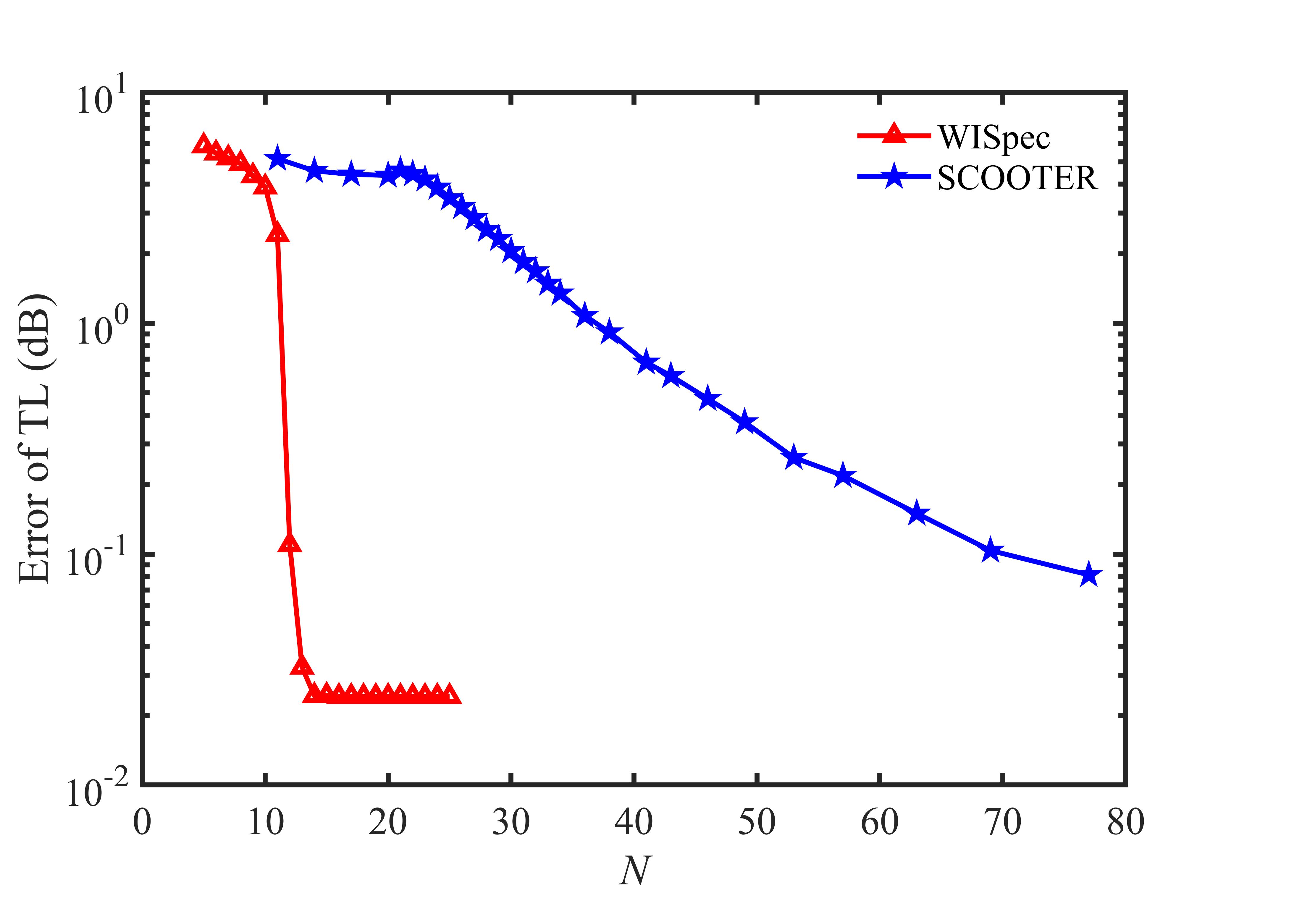}}
	\caption{Wavenumber spectrum of the pseudolinear-speed waveguide with a perfectly rigid bottom calculated by WISpec (a); wavenumber spectrum at a depth of $z=50$ m (b); sound fields calculated via the analytical solution (c) and WISpec (d); TLs along the $r$-direction at a depth of $z=50$ m (e); and the variation trend of the error in the TL field with $N$ (f).}
	\label{Figure4}
\end{figure}

A pseudolinear-speed waveguide is a waveguide whose sound speed profile has the following form \cite{Jensen2011}:
\begin{equation}
	c(z)=\sqrt{\frac{1}{az+b}}
\end{equation}
A pseudolinear-speed waveguide has an analytical solution involving Airy functions ($Ai(\cdot)$ and $Bi(\cdot)$) and their first derivatives ($Ai'(\cdot)$ and $Bi'(\cdot)$), and the horizontal wavenumbers $k_r$ are the roots of a transcendental equation \cite{Chowdhury2018}. In this example, the seabed is perfectly rigid. We take a sea depth $H=100$ m, $a=5.94\times10^{-10}$ s$^2$/m$^3$, and $b=4.16\times10^{-7}$ s$^2$/m$^3$; the sound source frequency $f=50$ Hz, $M=4096$, and the spectral truncation order used in WISpec is $N=15$. Table \ref{tab2} lists the discrete modes calculated using the three methods. The results of WISpec (see Fig.~\ref{Figure4a} and \ref{Figure4b}) are very consistent with the analytical solution. The agreement of the sound fields in Fig.~\ref{Figure4c} to \ref{Figure4e} also illustrates the reliablity of WISpec. Fig.~\ref{Figure4f} shows the variation trends of the errors of WISpec and SCOOTER for calculating the sound field as functions of the spectral truncation order or the number of grid points, from which it can be intuitively seen that WISpec achieves faster convergence and lower error. For a pseudolinear-speed waveguide with a free bottom, the same conclusion can be drawn, which is not repeated here.

\begin{table}[htbp]
	\centering
	\caption{Discrete modes of the pseudolinear-speed waveguide (unit: m$^{-1}$).}
		\begin{tabular}{cccc}
			\hline
			Mode order &Analytical solution &WISpec ($N$=15) &SCOOTER ($N$=15)   \\
			\hline
			1  & 0.2130 &0.2130 &0.2131\\
			2  & 0.2044 &0.2045 &0.2044\\
			3  & 0.1943 &0.1943 &0.1944\\
			4  & 0.1785 &0.1785 &0.1785\\
			5  & 0.1549 &0.1549 &0.1547\\
			6  & 0.1188 &0.1188 &0.1189\\
			7  & 0.0477 &0.0478 &0.0478\\
			\hline
		\end{tabular}
	\label{tab2}
\end{table}

\subsection{Pekeris waveguide}
\begin{figure}[htbp]
	\centering
	\subfigure[]{\includegraphics[width=8cm]{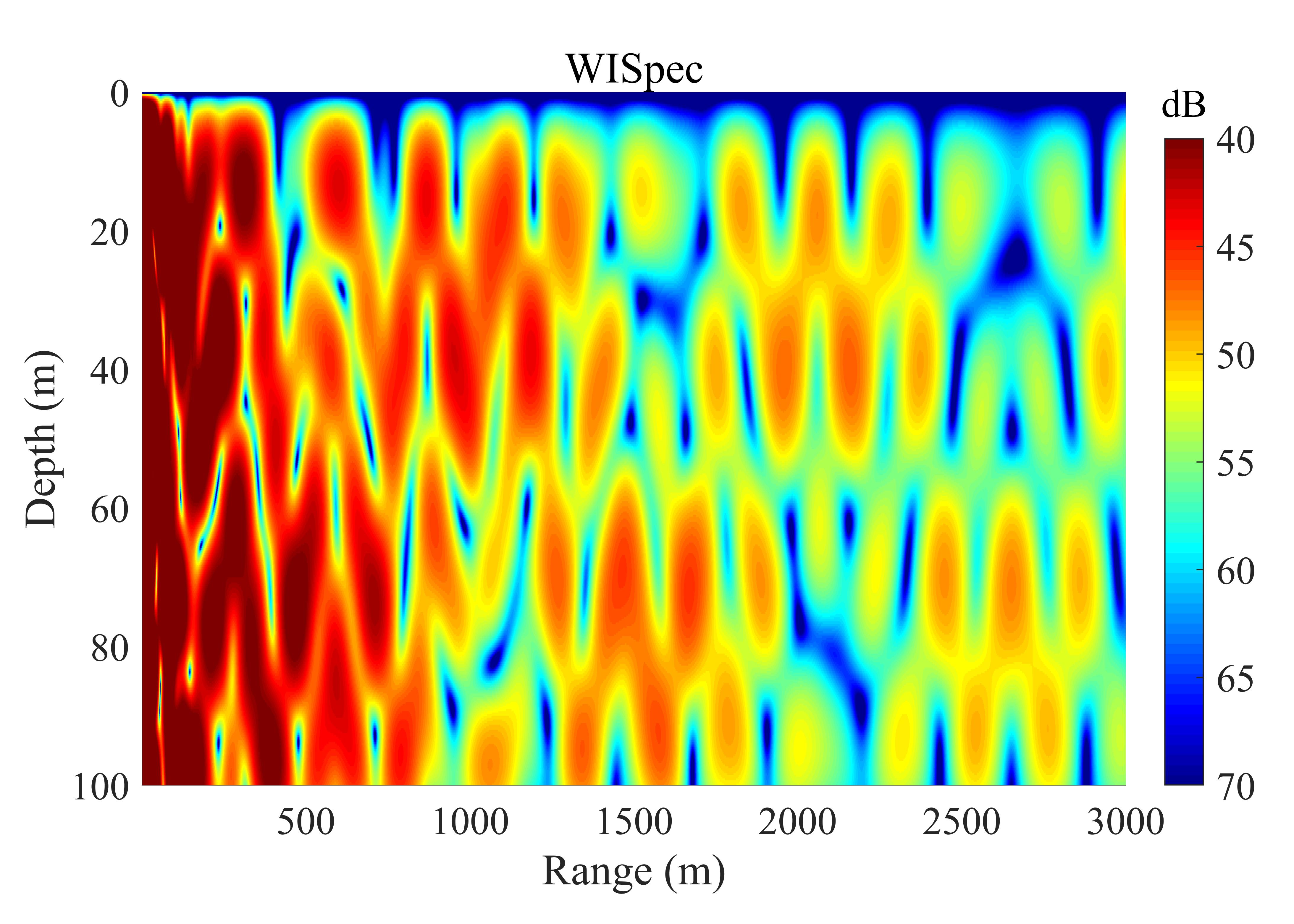}}
	\subfigure[]{\includegraphics[width=8cm]{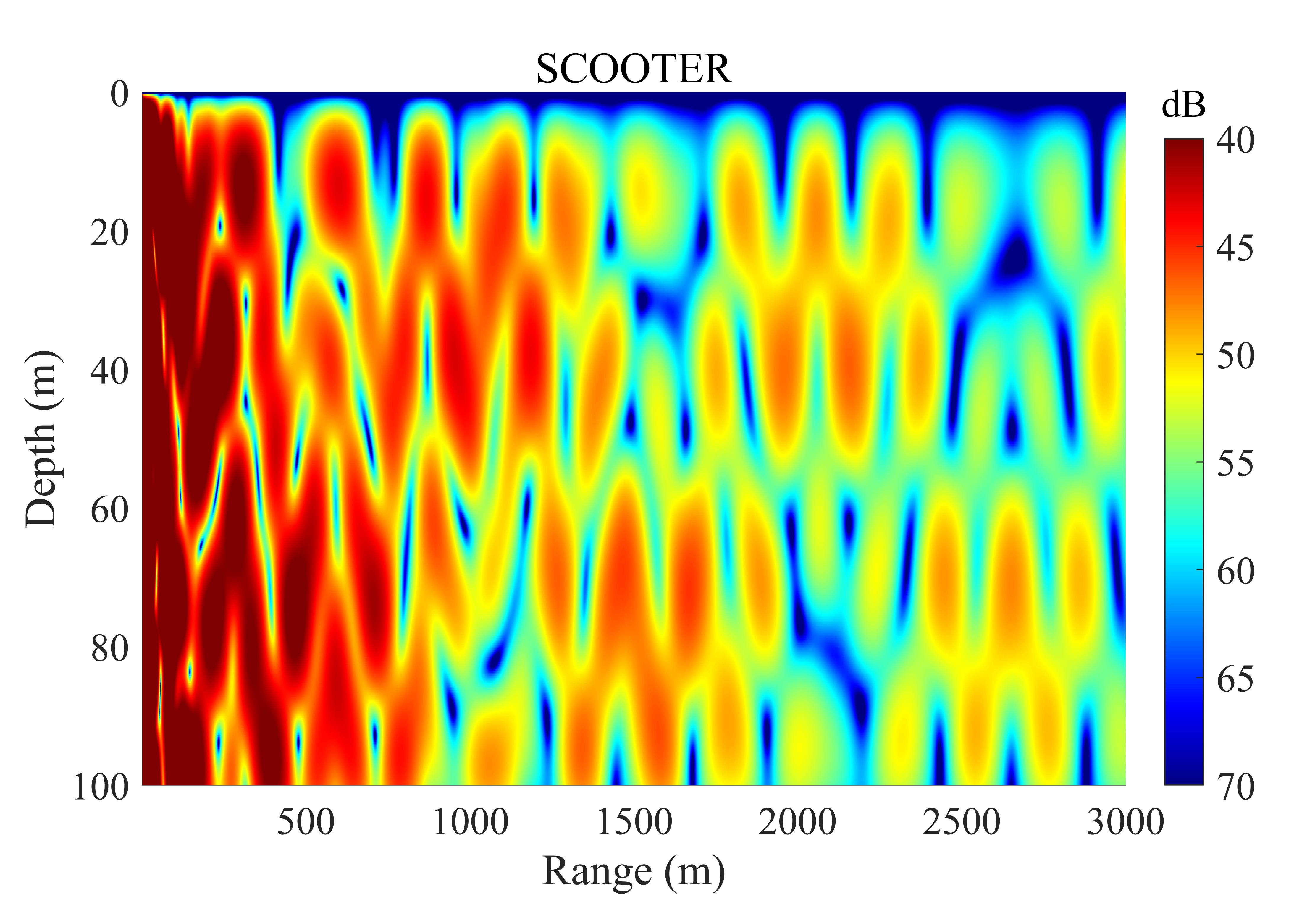}}\\
	\subfigure[]{\includegraphics[width=8cm]{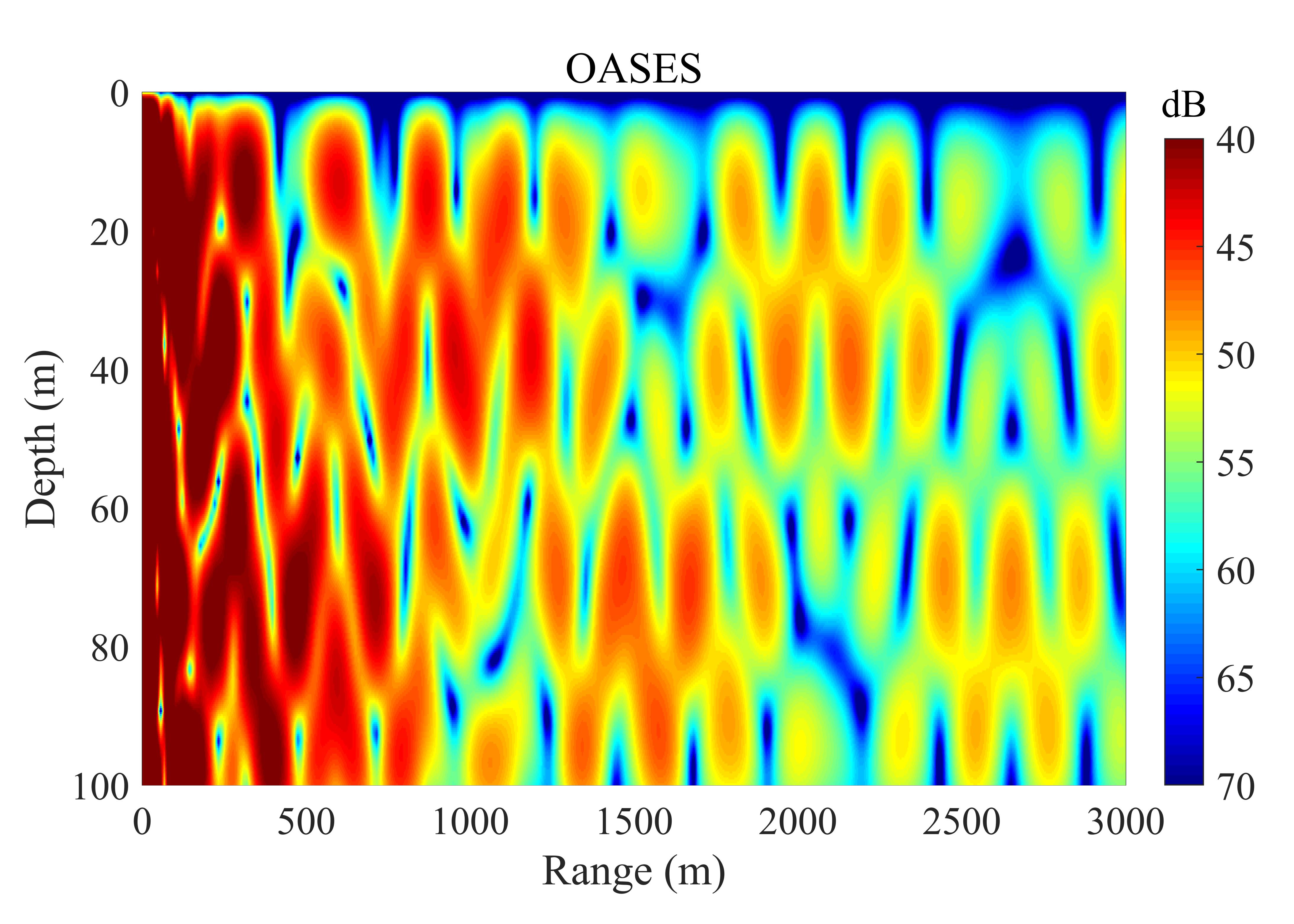}}
	\subfigure[]{\includegraphics[width=8cm]{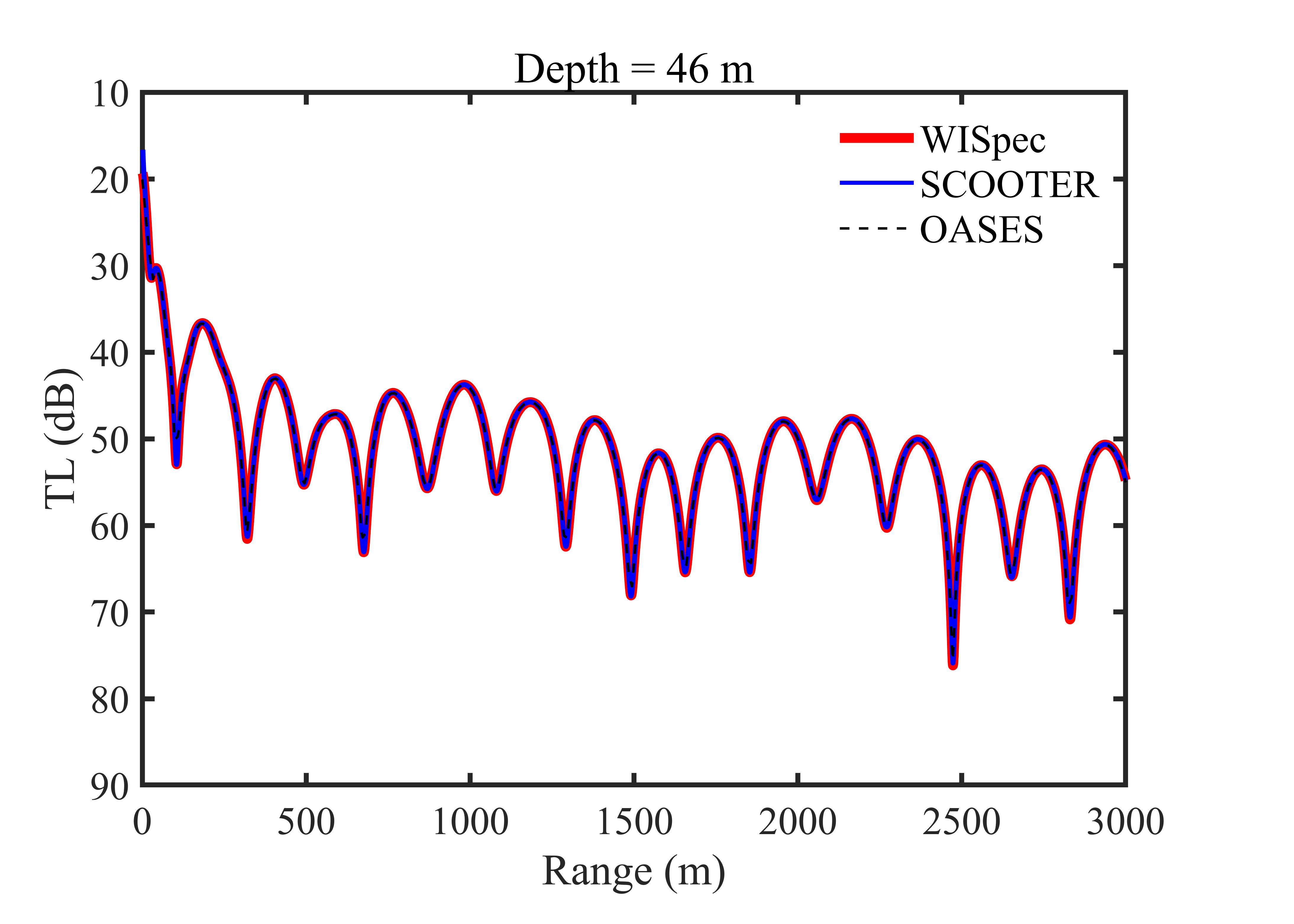}}
	\caption{Sound fields of the Pekeris waveguide calculated by WISpec (a), SCOOTER (b) and OASES (c); TLs along the $r$-direction at a depth of $z=46$ m (d).}
	\label{Figure5}
\end{figure}

The Pekeris waveguide is a classic waveguide in ocean acoustics; the ocean environment of the Pekeris waveguide consists of a layer of homogeneous water and an acoustic half-space below it. In this example, the same configuration as the ideal fluid waveguide is used, except that the sound source frequency is $f=50$ Hz, the density in the acoustic half-space is $\rho_\infty=1.5$ g/cm$^3$, the speed of sound is $c_\infty=2000$ m/s, and the attenuation is $\alpha_\infty=$ 0.5 dB/$\lambda$. We present the results of SCOOTER and OASES \cite{OASES} in Fig.~\ref{Figure5} for reference. OASES is a code to model seismic-acoustic propagation in horizontally stratified waveguides using wavenumber integration and the DGM approach. It is an upgraded version of SAFARI \cite{SAFARI}. The sound fields calculated by the three programs are basically the same. Whether on the sound fields or the TL-lines diagram, the three programs based on the wavenumber integration model have satisfactory consistency. This also proves the capability of WISpec to simulate waveguides with an acoustic half-space (including a continuous spectrum).

\begin{figure}[htbp]
	\centering
	\subfigure[]{\includegraphics[width=8cm]{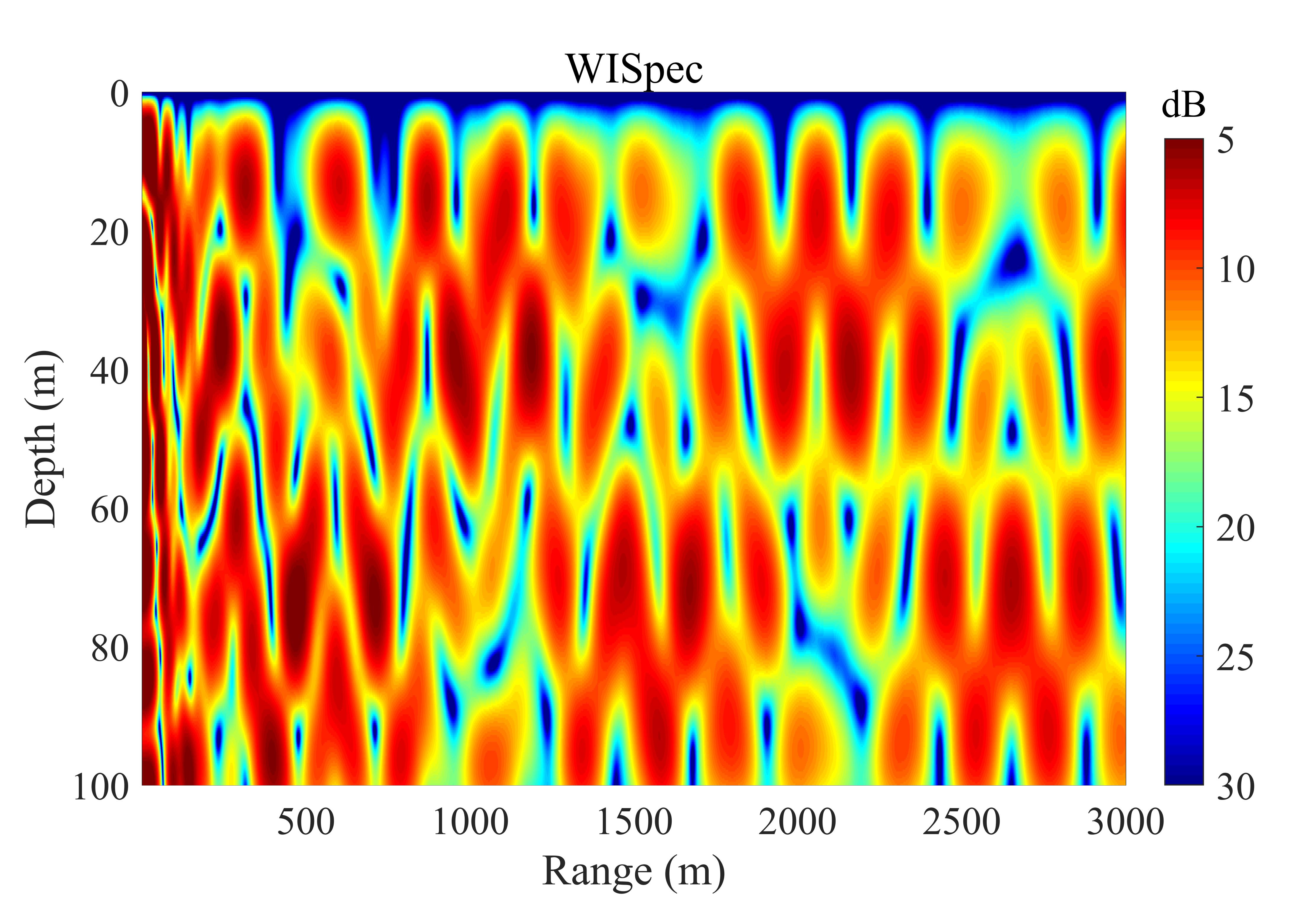}}
	\subfigure[]{\includegraphics[width=8cm]{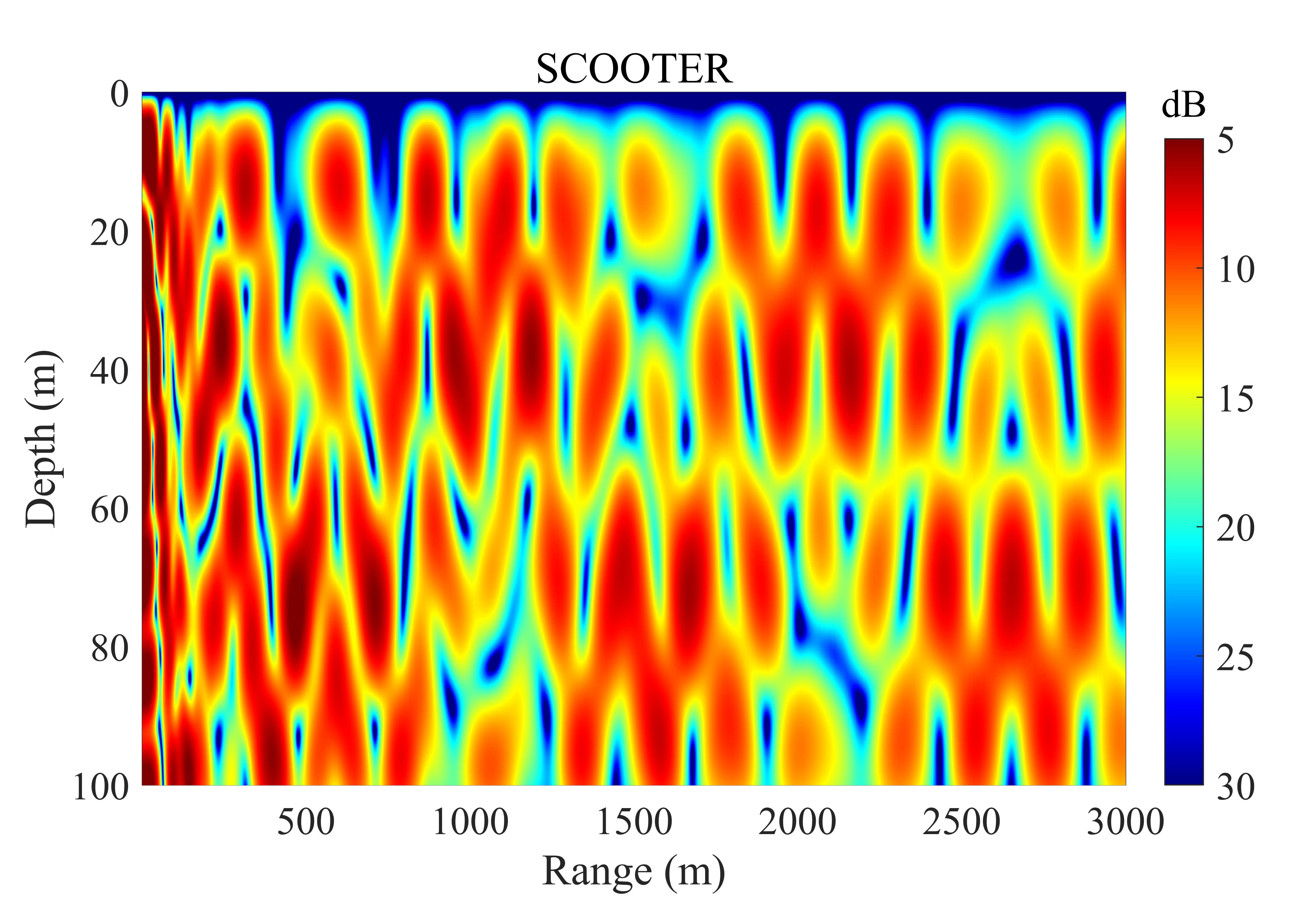}}\\
	\subfigure[]{\includegraphics[width=8cm]{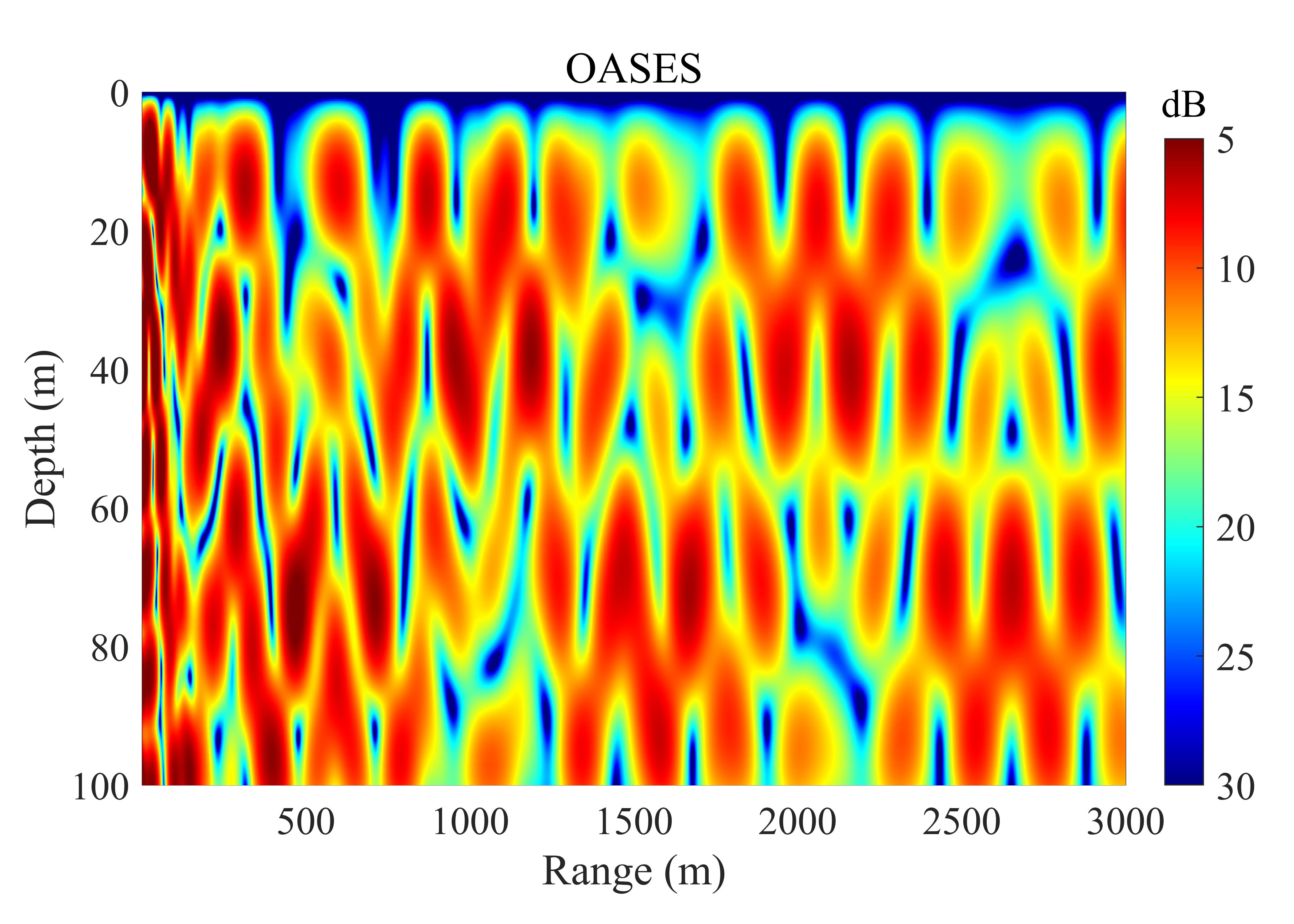}}
	\subfigure[]{\includegraphics[width=8cm]{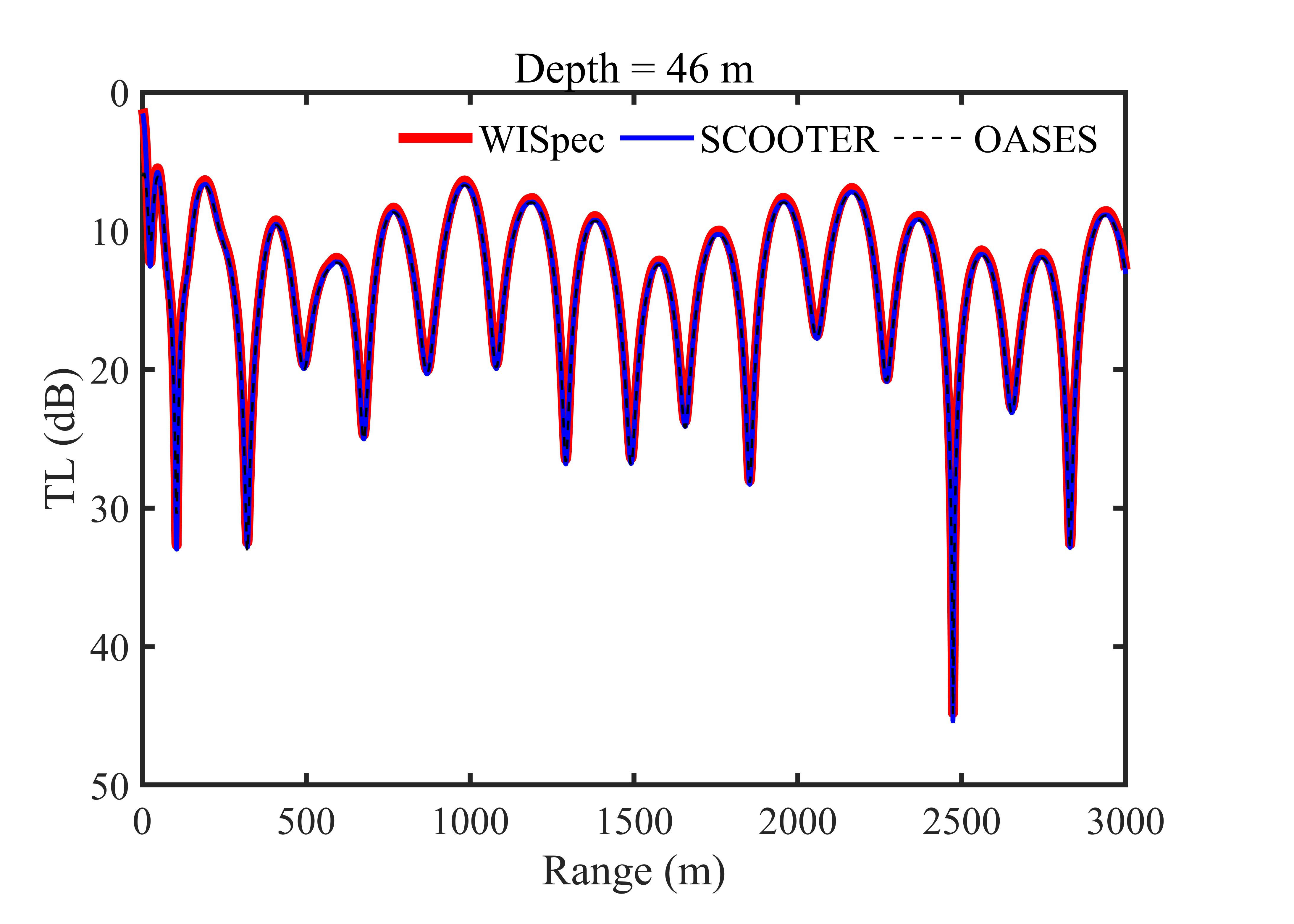}}
	\caption{Sound fields of the Pekeris waveguide of the line source calculated by WISpec (a), SCOOTER (b) and OASES (c); TLs along the $x$-direction at a depth of $z=46$ m (d).}
	\label{Figure6}
\end{figure}

In addition to point sources, WISpec can also calculate the sound field of line sources. Replacing the sound source with a line source in this example results in the sound field displayed in Fig.~\ref{Figure6}. The sound fields calculated by WISpec, SCOOTER and OASES are still very similar. To facilitate the comparison with SCOOTER and OASES, the sound field of the line source in WISpec is normalized as $p_{0}=\mathrm{i}\rho_\mathrm{s} \omega^{2}\mathcal{H}_{0}^{(1)}(1)/4$ instead of Eq.~\eqref{eq.18}.

\subsection{Bucker waveguide}
\begin{figure}[htbp]
	\centering
	\subfigure[]{\begin{tikzpicture}[node distance=2cm]
			\filldraw[very thick,cyan,opacity=0.4] (2,0)--(6.5,0)--(6.5,-5.25)--(2,-5.25)--cycle;
			\filldraw[very thick,brown,opacity=0.7] (2,-5.25)--(6.5,-5.25)--(6.5,-7)--(2,-7)--cycle;		
			\node at (1.8,0){$0$};
			\draw[very thick, ->](1.98,0)--(7,0) node[right]{$r$};
			\draw[very thick, ->](2,0)--(2,-7.5) node[below]{$z$};	
			\node at (1.7,-5.25){$H$};
			\filldraw [red] (2,-1) circle [radius=2.5pt];
			\draw[color=blue,very thick](4,0)--(3.5,-2.625);
			\draw[color=blue,very thick](3.5,-2.625)--(4,-5.25);
			\node at (3,-1){$\text{source}$};
			\node at (5.3,-0.5){$c=1500$ m/s};
			\node at (5.3,-2.625){$c=1498$ m/s};
			\node at (5.3,-5){$c=1500$ m/s};
			\node at (4,-6){$c_{\infty}=1505$ m/s};
			\node at (4,-6.5){$\rho_{\infty}=2.1$ g/cm$^3$};
			\node at (5.3,-1.5){$\rho=1$ g/cm$^3$};
		\end{tikzpicture}\label{Figure7a}}
	\subfigure[]{\begin{tikzpicture}[node distance=2cm,samples=1000,domain=0:5.25]
			\filldraw[very thick,cyan,opacity=0.4] (2,0)--(6.5,0)--(6.5,-5.25)--(2,-5.25)--cycle;
			\filldraw[very thick,cyan,opacity=0.7] (2,-5.25)--(6.5,-5.25)--(6.5,-7)--(2,-7)--cycle;		
			\node at (1.8,0){$0$};
			\draw[very thick, ->](1.98,0)--(7,0) node[right]{$r$};
			\draw[very thick, ->](2,0)--(2,-7.5) node[below]{$z$};	
			\node at (1.7,-5.25){$H$};
			\filldraw [red] (2,-0.5) circle [radius=2.5pt];
			\node at (2.6,-0.3){$\text{source}$};
			\node at (4,-6){$c_{\infty}=1600$ m/s};
			\node at (4,-6.5){$\rho_{\infty}=1$ g/cm$^3$};
			\node at (4,-2.5){$\rho=1$ g/cm$^3$};
			\draw[color=blue,very thick,smooth,rotate=-90] plot(\x,{22.5*(1+0.00737*((\x-1.6)/0.6-1+exp(-(\x-1.6)/0.6)))-20});
		\end{tikzpicture}\label{Figure7b}}
	\caption{Schematics of the ocean environment of the Bucker waveguide (a) and Munk waveguide (b).}
\end{figure}
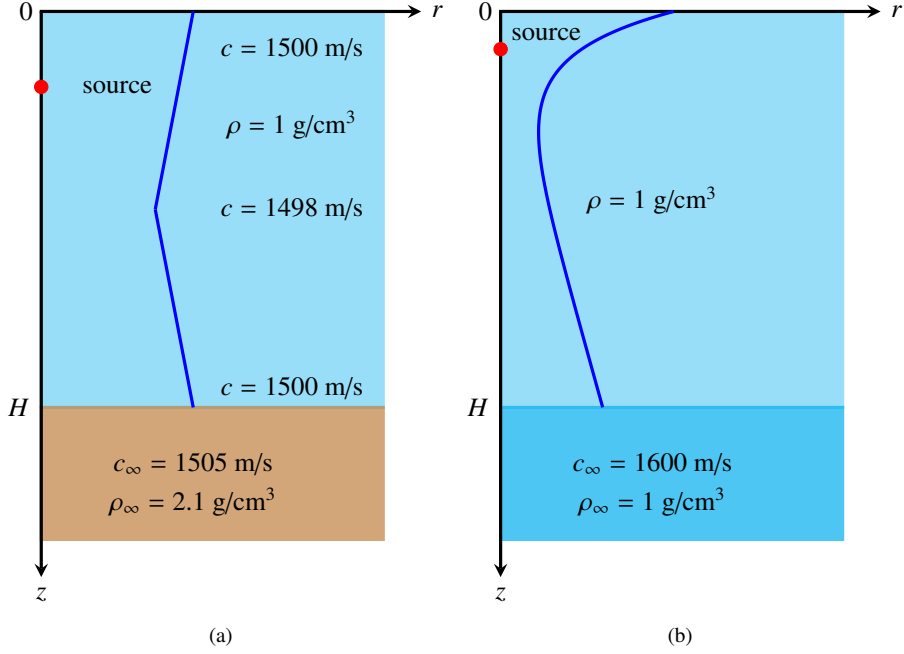

\begin{figure}[htbp]
	\centering
	\subfigure[]{\includegraphics[width=8cm]{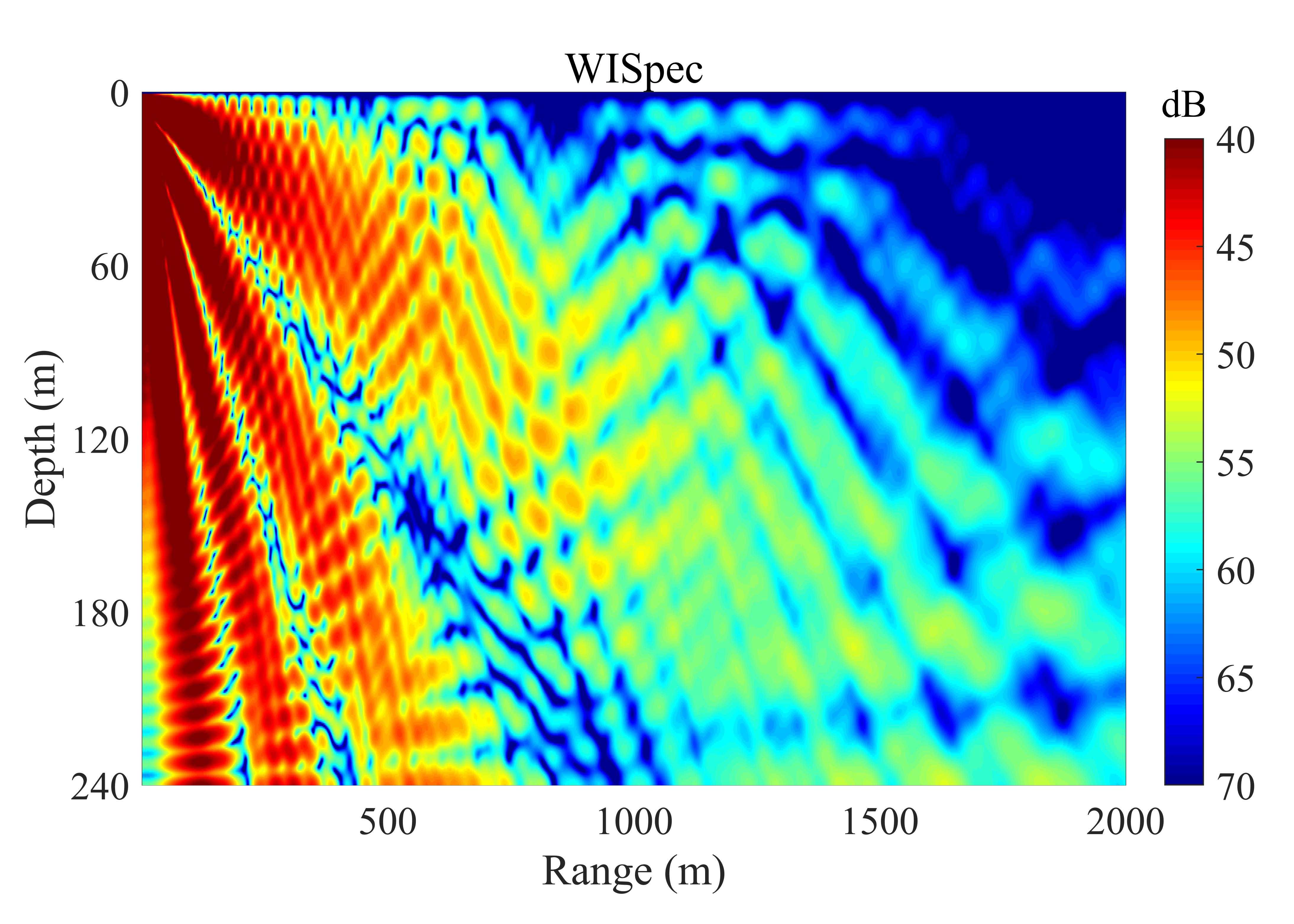}}
	\subfigure[]{\includegraphics[width=8cm]{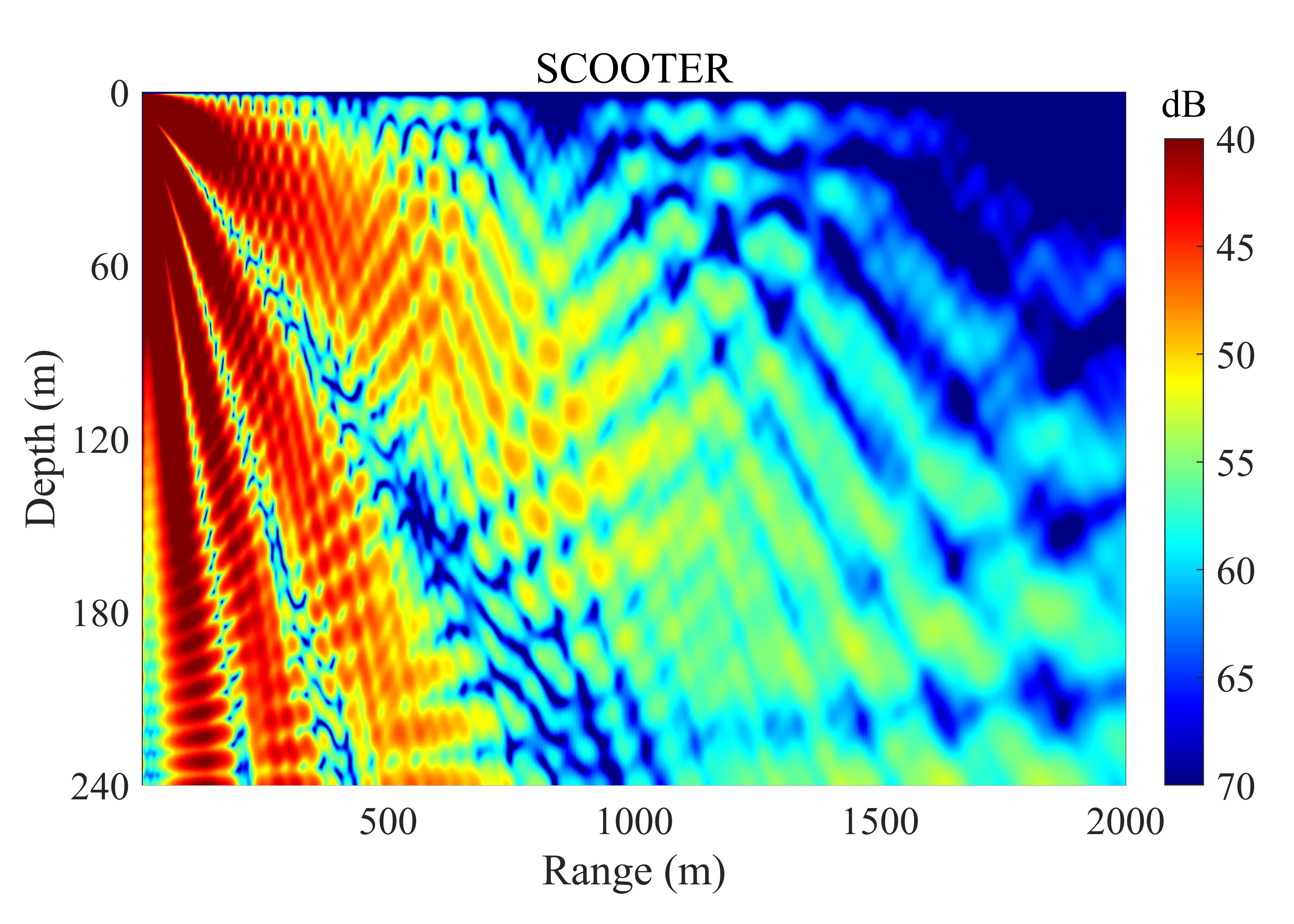}}\\
	\subfigure[]{\includegraphics[width=8cm]{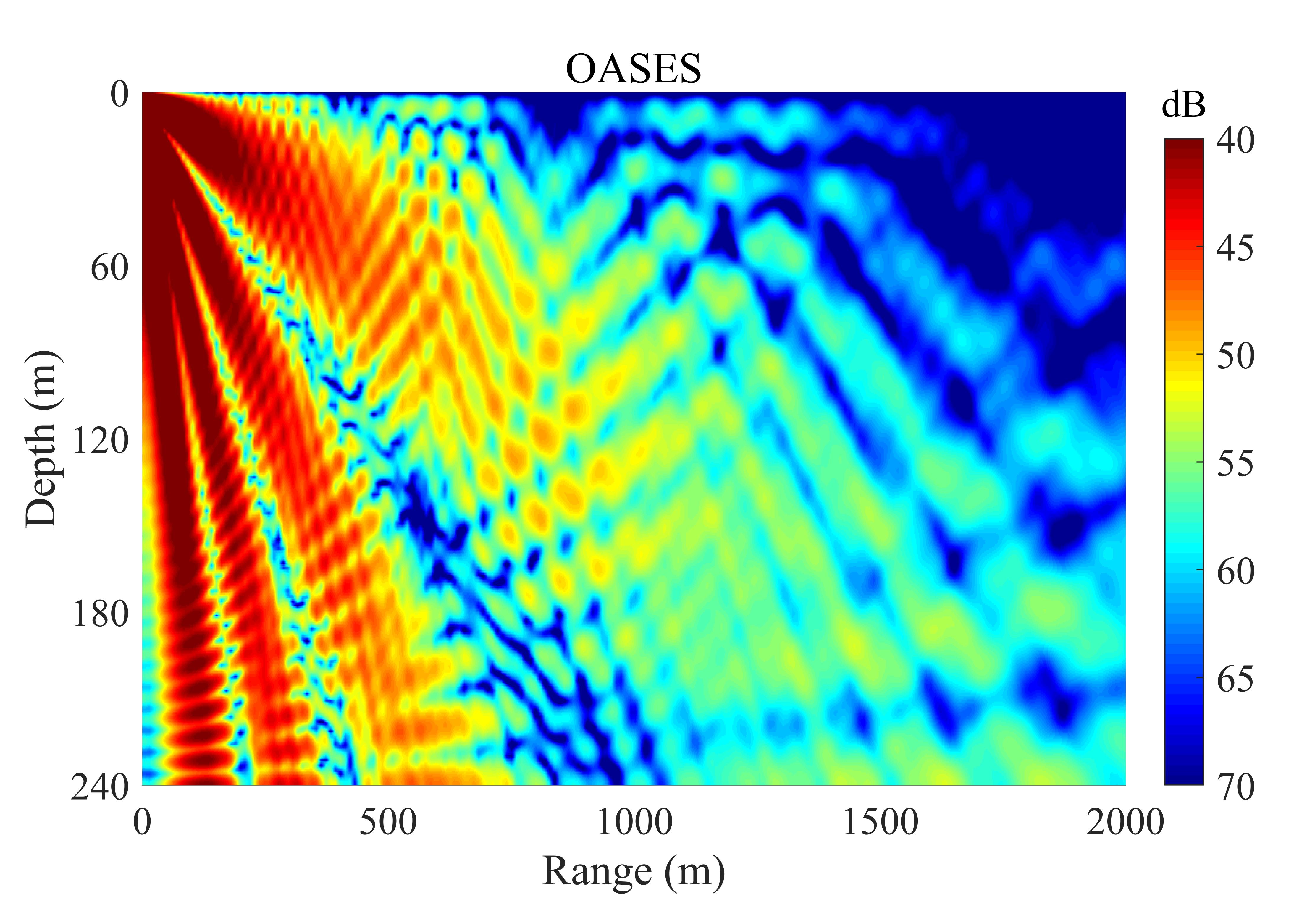}}
	\subfigure[]{\includegraphics[width=8cm]{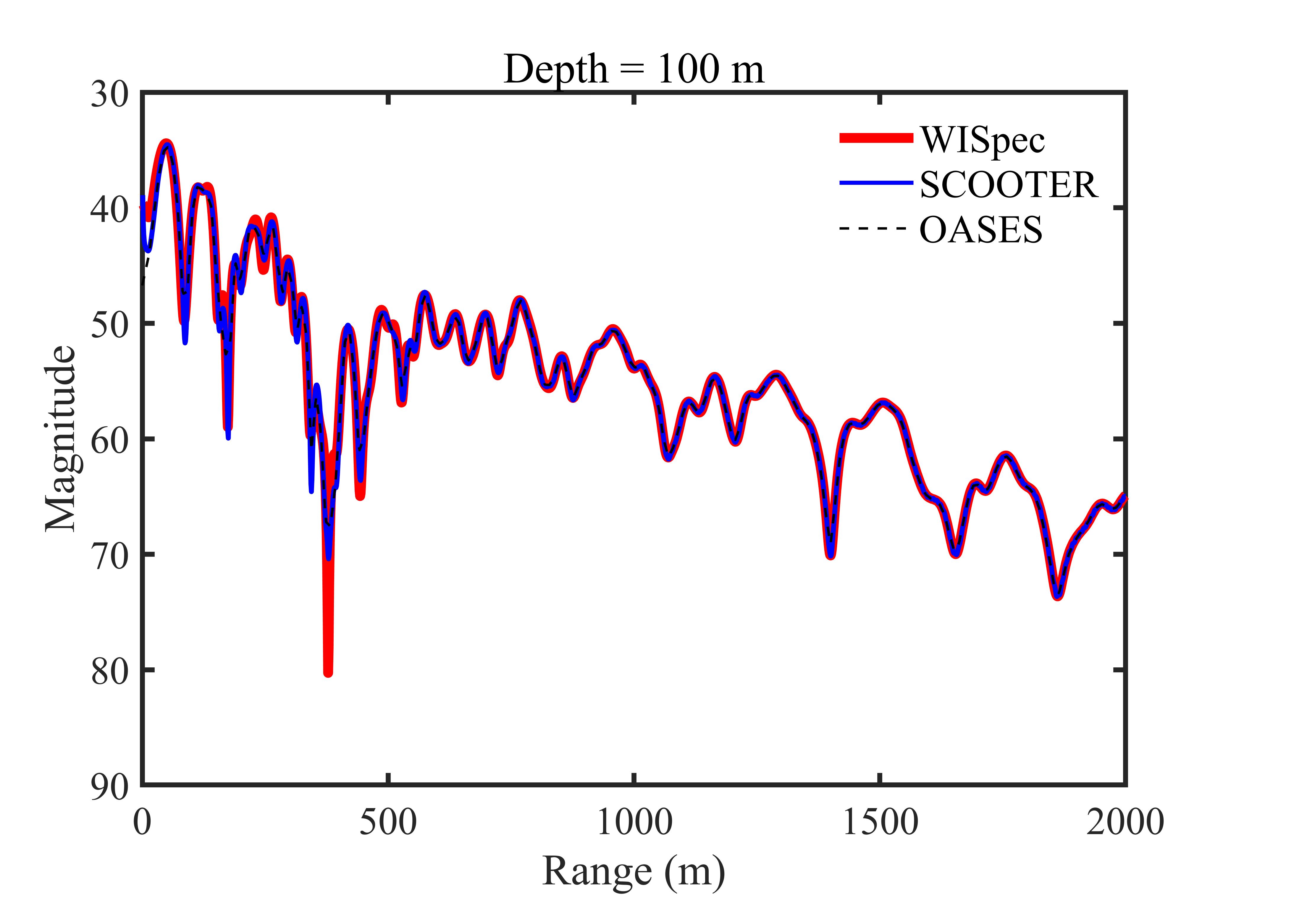}}
	\caption{Sound fields of the Bucker waveguide calculated by WISpec (a), SCOOTER (b) and OASES (c); TLs along the $r$-direction at a depth of $z=100$ m (d).}
	\label{Figure8}
\end{figure}

The Bucker waveguide is a benchmark for ocean acoustic propagation models \cite{Jensen2011}. As shown in Fig.~\ref{Figure7a}, the sound speed difference is very small, yielding a small number of propagating modes. However, this environment is characterized by a strong density difference at the sediment, which generates a large number of virtual modes close to the real wavenumber axis. Therefore, the normal mode model ignoring the continuous spectrum is not able to provide accurate results for the TL. However, wavenumber integration models are capable of providing an exact numerical prediction for the Bucker waveguide. In this experiment, the frequency of the sound source is taken as $f=100$ Hz, $z_\mathrm{s}=30$ m, and the sea depth is $H=240$ m. The number of discrete points in the wavenumber domain is taken as $M=4096$, the integral interval is $[0,2k_0]$, and the spectral truncation order is $N=40$. Fig.~\ref{Figure8} lists the sound fields of the Bucker waveguide calculated using WISpec, SCOOTER and OASES. The results of the three programs are in good agreement, except for the minor differences between the convergence area and the shadow area in the near field.

\subsection{Munk waveguide}
\begin{figure}[htbp]
	\centering
	\subfigure[]{\includegraphics[width=8cm]{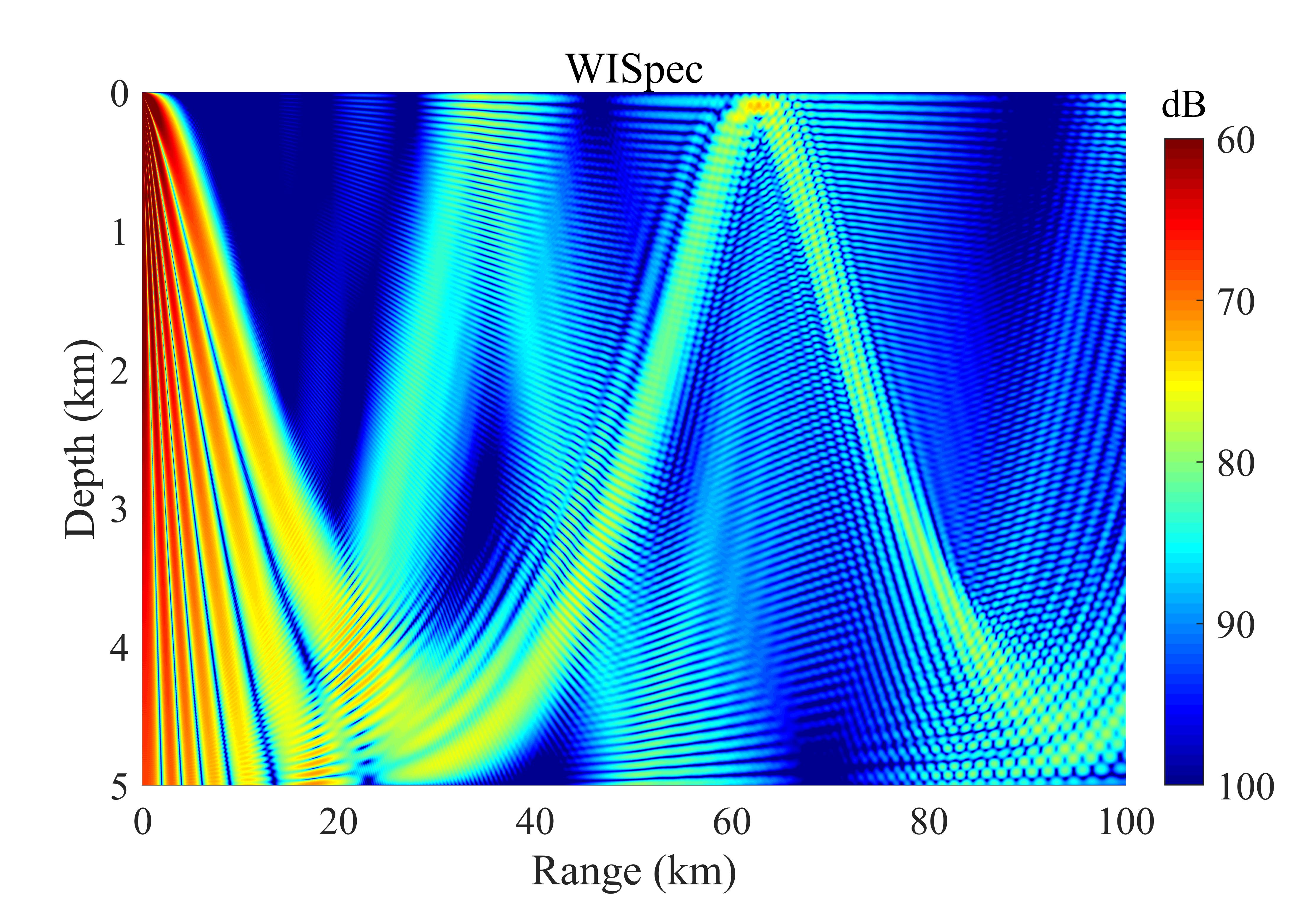}}
	\subfigure[]{\includegraphics[width=8cm]{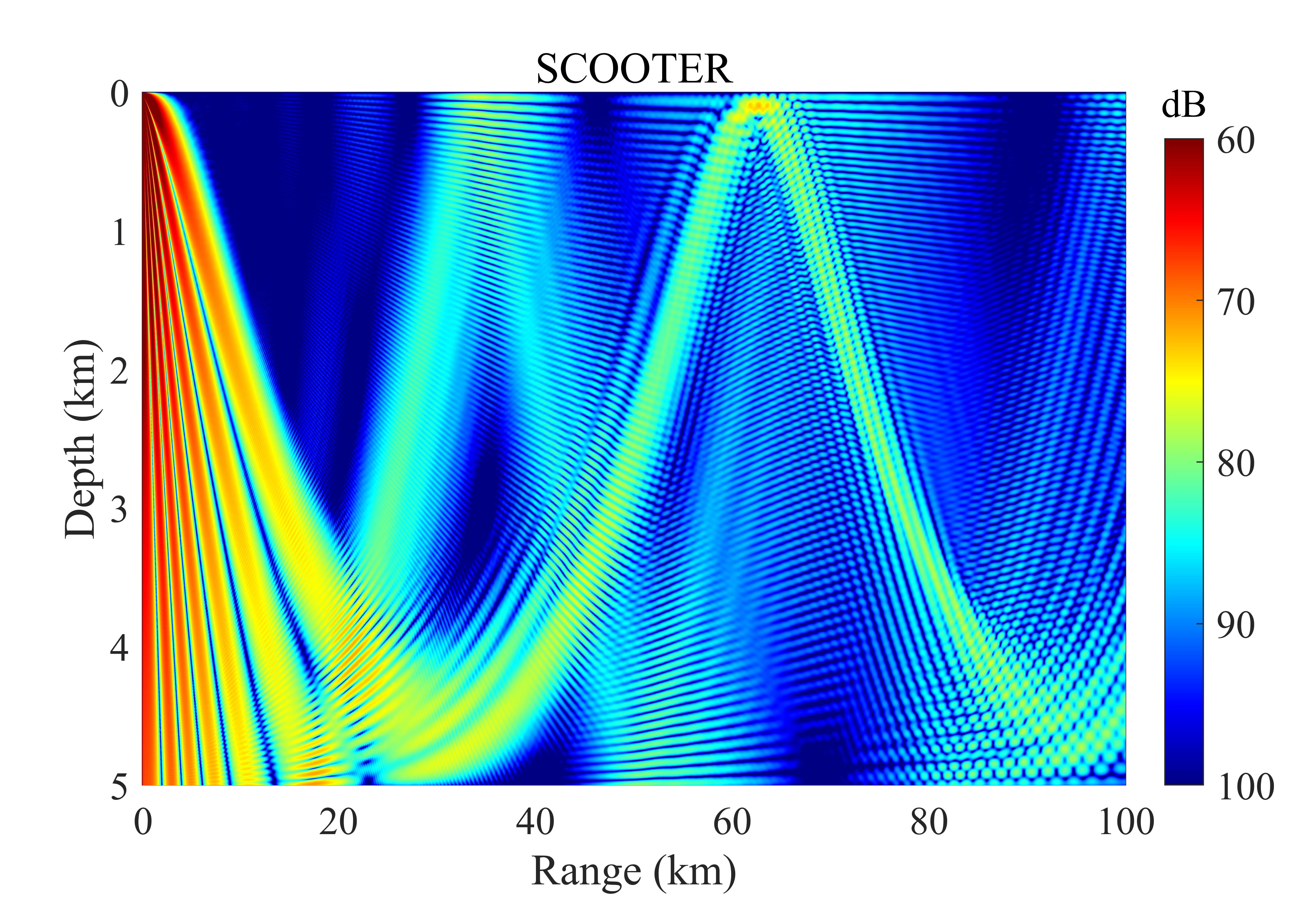}}\\
	\subfigure[]{\includegraphics[width=8cm]{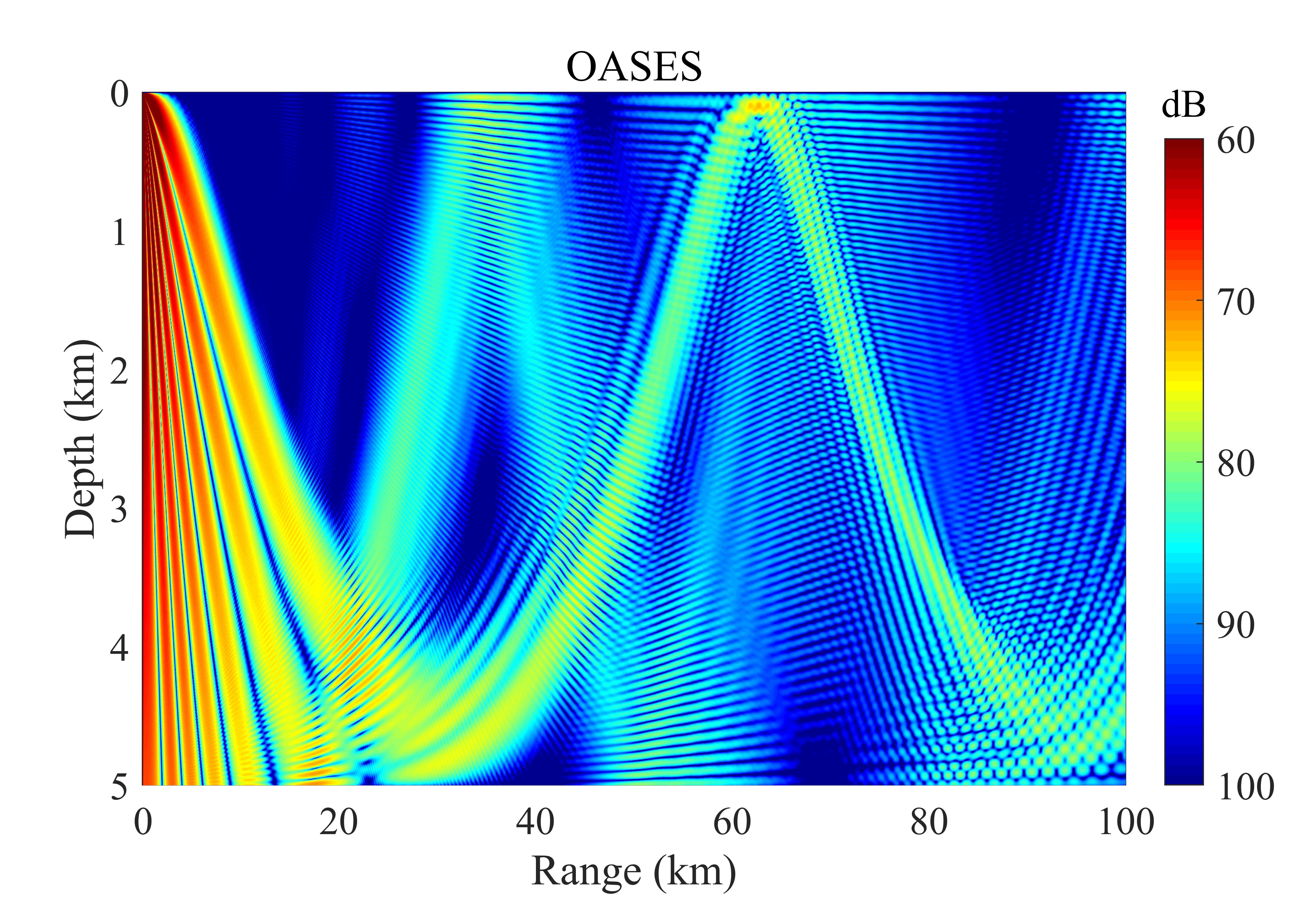}}
	\subfigure[]{\includegraphics[width=8cm]{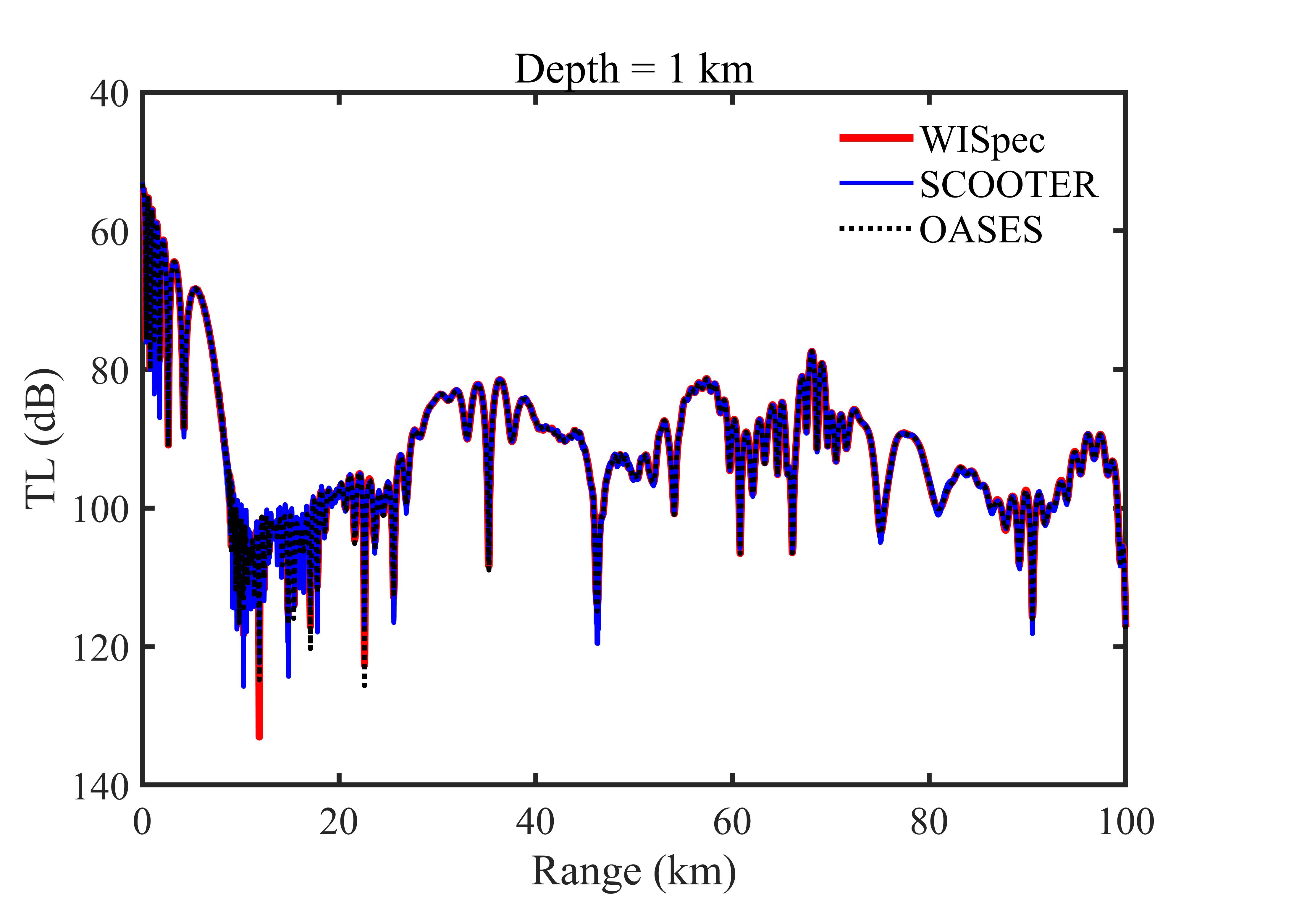}}
	\caption{Sound fields of the Munk waveguide calculated by WISpec (a), SCOOTER (b) and OASES (c); TLs along the $r$-direction at a depth of $z=1000$ m (d).}
	\label{Figure9}
\end{figure}

The Munk waveguide is a typical example of deep-sea acoustic propagation problems. Here, the ocean environment consists of a layer of seawater with a sound speed profile of the Munk profile and a homogeneous half-space below, as schematically shown in Figure~\ref{Figure7b}. In this experiment, the frequency of the sound source is $f=50$ Hz, $z_\mathrm{s}=100$ m, the sea depth is $H=5000$ m, and the sound speed profile is \cite{Jensen2011}:
\begin{equation}
	\begin{gathered}
		c(z)=1500.0[1.0+\epsilon(\tilde{z}-1+\exp{(-\tilde{z})})]\\
		\epsilon=0.00737,\quad \tilde{z}=(z-1300)/650.
	\end{gathered}
\end{equation}
The number of discrete points in the wavenumber domain is taken as $M=55000$, the integral interval is $[0,2k_0]$, and the spectral truncation order is $N=400$. Fig.~\ref{Figure9} illustrates the sound fields of the Munk waveguide calculated by WISpec, SCOOTER and OASES. The results of the three programs are very similar. From the TL-curves diagram, the results of WISpec and the other two software programs match very well. It is worth mentioning that, compared with WISpec and OASES, SCOOTER's result is slightly jittery at a range of 10-20 km.

To evaluate the computational speed of the algorithm proposed in this paper, we list the time consumption of solving the Green function in the above numerical experiments in Table \ref{tab3}. The test values for runtime are derived from the average of ten tests run per program on the same hardware platform (HUAWEI MateBookX Pro 2018, Intel i7-8550U CPU and 8 GB RAM). In addition, the times listed in the table are measured when the simulation results are reliable and the accuracy is approximately comparable. As seen from the table, WISpec has an approximately comparable runtime to SCOOTER and OASES for low frequencies and shallow seas and can achieve slightly faster speeds. However, for the deep-sea and high-frequency waveguide, the computational speed is significantly slower than for the above two mature programs. This also shows that WISpec has potential for further optimization.

\begin{table}[htbp]
	\centering
	\caption{Runtimes when solving the depth-dependent Green function in the numerical experiments (unit: s).}
		\begin{tabular}{lcccc}
			\hline
			{Case}&{Special configuration}&
			{WISpec}&
			{SCOOTER}&
			{OASES}\\
			\hline
			\multirow{2}{*}{Ideal fluid} &free seabed&0.513  & 0.583   &0.694  \\
			&rigid seabed&0.484  & 0.653   &0.703  \\
			\hline
			\multirow{2}{*}{Pseudolinear} &free seabed&0.814  &0.905   &1.001\\
			&rigid seabed&0.845  &0.916    &0.994\\
			\hline
			\multirow{2}{*}{Pekeris}
			&--         &0.818  & 0.899  & 0.746 \\
			&line source&0.819  & 0.903  & 0.726 \\
			\hline
			Bucker &--&1.649   & 1.164 &1.862\\
			\hline
			Munk &--&155.235  &70.882    &73.295  \\
			\hline
		\end{tabular}
	\label{tab3}
\end{table}

\section{Remarks and Summary}
\subsection{Remarks}
The above simulation experiments confirm that WISpec is a robust and accurate program and that the spectral method is effective in solving the depth-separated wave equation. From the above analysis, we can summarize the following features of the algorithm and program developed in this article:
\begin{enumerate}
	
	\item
	When applying the spectral method to solve the depth equation, as shown in Eq.~\eqref{eq.29}, there is no need to use piecewise linear elements to approximate the environmental parameters; i.e., the need to subdivide the environment into homogeneous sublayers is eliminated, thus avoiding error due to physical discretization in the vertical direction. The error of the spectral algorithm comes from the accuracy of spectral approximation. With the increase of truncation order, the error of spectral algorithm decreases exponentially. Therefore, the spectral algorithm can obtain accurate and reliable results.
	
	\item
	The discretization of the depth-separated wave equation by means of the Chebyshev--Tau spectral method yields a block diagonal matrix (see Eq.~\eqref{eq.33}), and in many cases, the Chebyshev matrix is quasidiagonally dominant; this sparsity enable efficient solution. The boundary and interface conditions are reflected in the global matrix, and the Green function of each layer can be obtained by solving it once, so the numerical stability of the spectral algorithm is good.
	
	\item
	The advantage of spectral method is that the error decreases exponentially with the increase of $N$, so the WISpec only needs a small truncation order to achieve reliable accuracy, which makes the scale of discrete matrix small. Compared with the existing wavenumber integration models, WISpec has the advantages of computational speed in low-frequency and shallow sea waveguides.
	
\end{enumerate}

\subsection{Summary}
In this paper, a novel approach for solving the depth-separated wave equation of the wavenumber integration model is developed, which uses the Chebyshev--Tau spectral method to discretize the depth equation and derive a large block diagonal matrix. For each additional layer of medium, the discretized block diagonal matrix grows accordingly. The problem is subsequently discretized into a system of linear equations. Solving the system of linear equations yields spectral coefficients, which are transformed by inverse Chebyshev transforms to obtain depth-dependent Green functions. A complete wavenumber integration model (WISpec) based on the proposed new approach is implemented. First, the Helmholtz equation is transformed into the wavenumber domain by means of Hankel or Fourier transformation. Since horizontally stratified media are considered, the wavenumber kernel function satisfies the depth-separated wave equation. The algorithm first samples the wavenumbers in a preset interval $[k_{\min},k_{\max}]$ and solves the Green functions $\Psi(k_r,z)$ in parallel for the discrete wavenumbers obtained through sampling. After the wavenumber kernel function is obtained, the inverse Hankel/Fourier transform is applied to synthesize the sound-pressure field in the physical space.

As far as we know, this algorithm is the first to use a spectral method to solve the depth-separated wave equation. Spectral methods use the idea of function approximation to control accuracy and the idea of the weighted residual to discretize the equation. When the ocean environment parameters are sufficiently smooth, the Green function solution converges exponentially. The results of numerical simulations verify the accuracy and reliability of the approach and model. The robust and high-precision Chebyshev--Tau spectral method avoids the possible instability of traditional algorithms for the depth-separated wave equation.

In terms of its application scope, this model requires that the ocean environment be independent in the $r/x$ direction, which limits the practicality of WISpec to a certain extent. Therefore, developing a high-precision wavenumber integration model based on the spectral method to solve the range-dependent waveguide has a bright future. In addition, elastic sediment is a more accurate model of the real ocean environment. In the future, WISpec can be further improved to enable the prediction of more complicated ocean acoustic fields.
	
\section*{Acknowledgments}
	
The authors thank Amlan Datta Chowdhury of the Indian Institute of Technology Madras for his valuable guidance on the analytical solution of the sound field for the waveguide with a pseudolinear sound speed profile. The authors also thank Zikai Gao of the National University of Defense Technology for his valuable guidance on the installation of OASES.
	
This work was supported by the National Natural Science Foundation of China under Grant 61972406 and by the National Key Research and Development Program of China under Grant 2016YFC1401800.

	\bibliographystyle{elsarticle-num}
	%\bibliography{refs.bib}	

\begin{thebibliography}{10}
\expandafter\ifx\csname url\endcsname\relax
  \def\url#1{\texttt{#1}}\fi
\expandafter\ifx\csname urlprefix\endcsname\relax\def\urlprefix{URL }\fi
\expandafter\ifx\csname href\endcsname\relax
  \def\href#1#2{#2} \def\path#1{#1}\fi

\bibitem{Jensen2011}
F.~B. Jensen, W.~A. Kuperman, M.~B. Porter, H.~Schmidt, Computational Ocean
  Acoustics, Springer-Verlag, New York, 2011.
\newblock \href {https://doi.org/10.1007/978-1-4419-8678-8}
  {\path{doi:10.1007/978-1-4419-8678-8}}.

\bibitem{Etter2018}
P.~C. Etter, Underwater Acoustic Modeling and Simulation, CRC Press, Boca
  Raton, USA, 2018.
\newblock \href {https://doi.org/10.1201/9781315166346}
  {\path{doi:10.1201/9781315166346}}.

\bibitem{Pekeris1948}
C.~L. Pekeris, Theory of propagation of explosive sound in shallow water,
  Geological Society of America Memoirs 27~(1) (1948) 1--117.
\newblock \href {https://doi.org/10.1130/mem27-2-p1}
  {\path{doi:10.1130/mem27-2-p1}}.

\bibitem{EJP1957}
W.~M. Ewing, W.~S. Jardetzky, F.~Press, Elastic Wave in Layered Media,
  McGraw-Hill, New York, 1957.
\newblock \href {https://doi.org/10.1785/BSSA0470030290}
  {\path{doi:10.1785/BSSA0470030290}}.

\bibitem{Thomson1950}
W.~T. Thomson, Transmission of elastic waves through a stratified solid medium,
  Journal of Applied Physics 21~(2) (1950) 89--93.
\newblock \href {https://doi.org/10.1063/1.1699629}
  {\path{doi:10.1063/1.1699629}}.

\bibitem{Haskell1953}
N.~A. Haskell, The dispersion of surface waves in multilayered media, Bulletin
  of the Seismological Society of America 43~(1) (1953) 17--34.
\newblock \href {https://doi.org/10.1785/BSSA0430010017}
  {\path{doi:10.1785/BSSA0430010017}}.

\bibitem{Kennett2009}
B.~L.~N. Kennett, Seismic Wave Propagation in Stratified Media, ANU Press,
  Canberra, 2009.
\newblock \href {https://doi.org/10.22459/SWPSM.05.2009}
  {\path{doi:10.22459/SWPSM.05.2009}}.

\bibitem{Kennett1974}
B.~L.~N. Kennett, Reflections, rays and reverberations, Bulletin of the
  Seismological Society of America 64~(6) (1974) 1685--1696.
\newblock \href {https://doi.org/10.1785/BSSA0640061685}
  {\path{doi:10.1785/BSSA0640061685}}.

\bibitem{Kennett1979}
B.~L.~N. Kennett, N.~J. Kerry, Seismic waves in a stratified half-space,
  Geophysical Journal of the Royal Astronomical Society 57~(3) (1979) 557--583.
\newblock \href {https://doi.org/10.1111/j.1365-246X.1979.tb06779.x}
  {\path{doi:10.1111/j.1365-246X.1979.tb06779.x}}.

\bibitem{Schmidt1985}
H.~Schmidt, F.~B. Jensen, A full wave solution for propagation in multilayered
  viscoelastic media with application to gaussian beam reflection at
  fluid–solid interfaces, The Journal of the Acoustical Society of America
  77~(3) (1985) 813--825.
\newblock \href {https://doi.org/10.1121/1.392050}
  {\path{doi:10.1121/1.392050}}.

\bibitem{SAFARI}
H.~Schmidt,
  \href{https://openlibrary.cmre.nato.int/handle/20.500.12489/281}{User's guide
  of {SAFARI}, seismo-acoustic fast field algorithm for range-independent
  environments} (1988).
\newline\urlprefix\url{https://openlibrary.cmre.nato.int/handle/20.500.12489/281}

\bibitem{OASES}
H.~Schmidt,
  \href{https://oalib-acoustics.org/models-and-software/wavenumber-integration/}{{OASES}:
  User guide and reference manual (version 3.1)} (2020).
\newline\urlprefix\url{https://oalib-acoustics.org/models-and-software/wavenumber-integration/}

\bibitem{Schmidt1986}
H.~Schmidt, G.~J. Tango, Efficient global matrix approach to the computation of
  synthetic seismograms, Geophysical Journal International 84~(2) (1986)
  331--359.
\newblock \href {https://doi.org/10.1111/j.1365-246X.1986.tb04359.x}
  {\path{doi:10.1111/j.1365-246X.1986.tb04359.x}}.

\bibitem{Jensen1998}
F.~B. Jensen, On the use of stair steps to approximate bathymetry changes in
  ocean acoustic models, The Journal of the Acoustical Society of America
  104~(3) (1998) 1310--1315.
\newblock \href {https://doi.org/10.1121/1.424340}
  {\path{doi:10.1121/1.424340}}.

\bibitem{Orszag1972}
S.~A. Orszag, Comparison of pseudospectral and spectral approximation, Studies
  in Applied Mathematics L1~(3) (1972) 253--259.
\newblock \href {https://doi.org/10.1002/sapm1972513253}
  {\path{doi:10.1002/sapm1972513253}}.

\bibitem{Gottlieb1977}
D.~Gottlieb, S.~A. Orszag, Numerical Analysis of Spectral Methods, Theory and
  Applications, Society for Industrial and Applied Mathematics, Philadelphia,
  USA, 1977.
\newblock \href {https://doi.org/10.1137/1.9781611970425}
  {\path{doi:10.1137/1.9781611970425}}.

\bibitem{Canuto1988}
C.~Canuto, M.~Y. Hussaini, A.~Quarteroni, T.~A. Zang, Spectral Methods in Fluid
  Dynamics, Spring-Verlag, Berlin, Germany, 1988.
\newblock \href {https://doi.org/10.1007/978-3-642-84108-8}
  {\path{doi:10.1007/978-3-642-84108-8}}.

\bibitem{Guoby1998}
B.~Guo, Spectral Methods and Their Applications, World Scientific, 1998.
\newblock \href {https://doi.org/10.1142/3662} {\path{doi:10.1142/3662}}.

\bibitem{Boyd2001}
J.~P. Boyd, {Chebyshev} and {Fourier} Spectral Methods, Second Edition, Dover,
  New York, USA, 2001.

\bibitem{Canuto2006}
C.~Canuto, M.~Y. Hussaini, A.~Quarteroni, T.~A. Zang, Spectral Methods
  Fundamentals in Single Domains, Spring-Verlag, Berlin, German, 2006.
\newblock \href {https://doi.org/10.1007/978-3-540-30726-6}
  {\path{doi:10.1007/978-3-540-30726-6}}.

\bibitem{Jshen2011}
J.~Shen., T.~Tang., L.~Wang., Spectral Methods Algorithms, Analysis and
  Applications, Springer-Verlag, Berlin, German, 2011.
\newblock \href {https://doi.org/10.1007/978-3-540-71041-7}
  {\path{doi:10.1007/978-3-540-71041-7}}.

\bibitem{Sabatini2019}
R.~Sabatini, P.~Cristini, A multi-domain collocation method for the accurate
  computation of normal modes in open oceanic and atmospheric waveguides, Acta
  Acustica United with Acustica 105 (2019) 464--474.
\newblock \href {https://doi.org/10.3813/AAA.919328}
  {\path{doi:10.3813/AAA.919328}}.

\bibitem{Tuhw2021c}
Y.~Wang, H.~Tu, W.~Liu, W.~Xiao, Q.~Lan, Two {Chebyshev} spectral methods for
  solving normal modes in atmospheric acoustics, Entropy 23 (2021) 705.
\newblock \href {https://doi.org/10.3390/e23060705}
  {\path{doi:10.3390/e23060705}}.

\bibitem{Dzieciuch1993}
M.~A. Dzieciuch, Numerical solution of the acoustic wave equation using
  chebyshev polynomials with application to global acoustics, in: Proceedings
  of Oceans, IEEE, Victoria, BC, Canada, October 18--21, 1993, pp. 267--271.
\newblock \href {https://doi.org/10.1109/OCEANS.1993.326000}
  {\path{doi:10.1109/OCEANS.1993.326000}}.

\bibitem{Tuhw2020a}
H.~Tu, Y.~Wang, W.~Liu, X.~Ma, W.~Xiao, Q.~Lan, A {Chebyshev} spectral method
  for normal mode and parabolic equation models in underwater acoustics,
  Mathematical Problems in Engineering (2020) 7461314\href
  {https://doi.org/10.1155/2020/7461314} {\path{doi:10.1155/2020/7461314}}.

\bibitem{Tuhw2020b}
H.~Tu, Y.~Wang, Q.~Lan, W.~Liu, W.~Xiao, S.~Ma, A {Chebyshev--Tau} spectral
  method for normal modes of underwater sound propagation with a layered marine
  environment, Journal of Sound and Vibration 492 (2021) 115784.
\newblock \href {https://doi.org/10.1016/j.jsv.2020.115784}
  {\path{doi:10.1016/j.jsv.2020.115784}}.

\bibitem{Tuhw2021d}
H.~Tu, Y.~Wang, Q.~Lan, W.~Liu, W.~Xiao, S.~Ma, Applying a {Legendre}
  collocation method based on domain decomposition to calculate underwater
  sound propagation in a horizontally stratified environment, Journal of Sound
  and Vibration 511 (2021) 116364.
\newblock \href {https://doi.org/10.1016/j.jsv.2021.116364}
  {\path{doi:10.1016/j.jsv.2021.116364}}.

\bibitem{NM-CT}
H.~Tu,
  \href{https://oalib-acoustics.org/models-and-software/normal-modes/}{{NM-CT}:
  A {Chebyshev--Tau} spectral method for normal modes of underwater sound
  propagation with a layered marine environment in {Matlab} and {Fortran}}
  (2020).
\newline\urlprefix\url{https://oalib-acoustics.org/models-and-software/normal-modes/}

\bibitem{rimLG}
R.~B. Evans,
  \href{https://oalib-acoustics.org/models-and-software/normal-modes/}{{rimLG}:
  A {Legendre--Galerkin} technique for differential eigenvalue problems with
  complex and discontinuous coefficients, arising in underwater acoustics}
  (2020).
\newline\urlprefix\url{https://oalib-acoustics.org/models-and-software/normal-modes/}

\bibitem{MultiLC}
H.~Tu,
  \href{https://oalib-acoustics.org/models-and-software/normal-modes/}{{MultiLC}:
  A {Legendre} collocation method based on domain decomposition to calculate
  underwater sound propagation in a horizontally stratified environment in
  {Matlab} and {Fortran}} (2021).
\newline\urlprefix\url{https://oalib-acoustics.org/models-and-software/normal-modes/}

\bibitem{Tuhw2022a}
H.~Tu, Y.~Wang, C.~Yang, W.~Liu, W.~Xiao, A {Chebyshev--Tau} spectral method
  for coupled modes of underwater sound propagation in range-dependent ocean
  environments, arXiv.org (2022).
\newblock \href {https://doi.org/arXiv:2111.09493}
  {\path{doi:arXiv:2111.09493}}.

\bibitem{Tuhw2022b}
H.~Tu, Y.~Wang, C.~Yang, X.~Wang, S.~Ma, W.~Xiao, W.~Liu, A novel algorithm to
  solve for an underwater line source sound field based on coupled modes and a
  spectral method, Journal of Computational Physics 468 (2022) 111478.
\newblock \href {https://doi.org/10.1016/j.jcp.2022.111478}
  {\path{doi:10.1016/j.jcp.2022.111478}}.

\bibitem{Tuhw2022c}
H.~Tu, Y.~Wang, W.~Liu, C.~Yang, J.~Qin, S.~Ma, X.~Wang, Application of a
  spectral method to simulate quasi-three-dimensional underwater acoustic
  fields, Journal of Sound and Vibration 545 (2023) 117421.
\newblock \href {https://doi.org/10.1016/j.jsv.2022.117421}
  {\path{doi:10.1016/j.jsv.2022.117421}}.

\bibitem{Tuhw2023a}
H.~Tu, Y.~Wang, Y.~Zhang, H.~Liao, W.~Liu, Parallel numerical simulation of
  weakly range-dependent ocean acoustic waveguides by adiabatic modes based on
  a spectral method, Physics of Fluids 35~(1) (2023) 017119.
\newblock \href {https://doi.org/10.1063/5.0131771}
  {\path{doi:10.1063/5.0131771}}.

\bibitem{Tuhw2021a}
Y.~Wang, H.~Tu, W.~Liu, W.~Xiao, Q.~Lan, Application of a {Chebyshev}
  collocation method to solve a parabolic equation model of underwater acoustic
  propagation, Acoustics Australia (2021) 1--12\href
  {https://doi.org/10.1007/s40857-021-00218-5}
  {\path{doi:10.1007/s40857-021-00218-5}}.

\bibitem{Tuhw2021b}
H.~Tu, Y.~Wang, X.~Ma, X.~Zhu, Applying the {Chebyshev--Tau} spectral method to
  solve the parabolic equation model of wide-angle rational approximation in
  ocean acoustics, Journal of Theoretical and Computational Acoustics 30~(2)
  (2022).
\newblock \href {https://doi.org/10.1142/S2591728521500134}
  {\path{doi:10.1142/S2591728521500134}}.

\bibitem{SMPE}
H.~Tu,
  \href{https://oalib-acoustics.org/models-and-software/parabolic-equation/}{{SMPE}:
  Two spectral methods for solving the range-independent parabolic equation
  model in ocean acoustics} (2021).
\newline\urlprefix\url{https://oalib-acoustics.org/models-and-software/parabolic-equation/}

\bibitem{Luowy2016}
W.~Luo, X.~Yu, X.~Yang, R.~Zhang, Analytical solution based on the wavenumber
  integration method for the acoustic field in a {Pekeris} waveguide, Chinese
  Physics B 25~(4) (2016) 044302.
\newblock \href {https://doi.org/10.1088/1674-1056/25/4/044302}
  {\path{doi:10.1088/1674-1056/25/4/044302}}.

\bibitem{Lanczos1938}
C.~Lanczos, Trigonometric interpolation of empirical and analytical functions,
  Journal of Mathematical Physics 17 (1938) 123--199.

\bibitem{Min2005}
M.~S. Min, D.~Gottlieb, Domain decomposition spectral approximations for an
  eigenvalue problem with a piecewise constant coefficient, SIAM Journal on
  Numerical Analysis 43 (2005) 502--520.
\newblock \href {https://doi.org/10.1137/s0036142903423836}
  {\path{doi:10.1137/s0036142903423836}}.

\bibitem{SCOOTER}
M.~B. Porter,
  \href{https://oalib-acoustics.org/models-and-software/wavenumber-integration/}{{SCOOTER}:
  A finite element {FFP} code} (2022).
\newline\urlprefix\url{https://oalib-acoustics.org/models-and-software/wavenumber-integration/}

\bibitem{Chowdhury2018}
A.~D. Chowdhury, C.~P. Vendhan, S.~K. Bhattacharyya, S.~Mudaliar, A
  {Rayleigh-Ritz} model for the depth eigenproblem of heterogeneous pekeris
  waveguides, ACTA ACUSTICA UNITED WITH ACUSTICA 104 (2018) 597--610.
\newblock \href {https://doi.org/10.3813/aaa.919200}
  {\path{doi:10.3813/aaa.919200}}.

\end{thebibliography}

\end{document}